\documentclass[pdflatex,sn-mathphys-num]{sn-jnl}

\usepackage{amsmath,amssymb,mathtools}
\usepackage{amsthm}
\usepackage{enumitem}
\usepackage{microtype}
\usepackage{graphicx}

\numberwithin{equation}{section}
\theoremstyle{plain}
\newtheorem{theorem}{Theorem}[section]
\newtheorem{proposition}[theorem]{Proposition}
\newtheorem{lemma}[theorem]{Lemma}
\newtheorem{corollary}[theorem]{Corollary}

\theoremstyle{definition}

\theoremstyle{remark}
\newtheorem{remark}[theorem]{Remark}

\newcommand{\1}{\mathbf 1}
\newcommand{\BV}{\mathrm{BV}}
\newcommand{\TV}{\mathrm{TV}}

\usepackage{xcolor}
\hypersetup{hidelinks}

\newcommand{\compactdisplays}{%
  \setlength{\abovedisplayskip}{4pt plus 1pt minus 1pt}%
  \setlength{\belowdisplayskip}{4pt plus 1pt minus 1pt}%
  \setlength{\abovedisplayshortskip}{0pt plus 1pt}%
  \setlength{\belowdisplayshortskip}{2pt plus 1pt minus 1pt}%
}
\AddToHook{cmd/normalsize/after}{\compactdisplays}
\AddToHook{cmd/small/after}{\compactdisplays}
\AddToHook{cmd/footnotesize/after}{\compactdisplays}
\AtBeginDocument{\compactdisplays}
\AtBeginDocument{\raggedbottom\setlength{\parskip}{0pt}\setlength{\lineskip}{0pt}}
\makeatletter
\def\thm@space@setup{\thm@preskip=4pt\thm@postskip=4pt}
\g@addto@macro\ps@headings{\let\@oddhead\@empty\let\@evenhead\@empty}
\g@addto@macro\ps@titlepage{\let\@oddhead\@empty\let\@evenhead\@empty}
\makeatother
\AtBeginDocument{\pagestyle{headings}}

\begin{document}

\title{Conditional limit theorems for open intermittent dynamical systems}
\author*[1]{\fnm{Can} \sur{Wang}}\email{12631038@mail.sustech.edu.cn}
\author*[1]{\fnm{Yaofeng} \sur{Su}}\email{suyaofeng@sustech.edu.cn}
\affil[1]{\orgdiv{Department of Mathematics}, \orgname{Southern University of Science and Technology}, \orgaddress{\city{Shenzhen}, \country{China}}}
\abstract{
		In this paper, we study conditional limit theorems for open intermittent systems endowed with big/small holes. The initial measure for the system can be a broad class of (non)invariant probability, including SRB and Lebesgue measures. Unlike the limit theorems for closed systems, we obtain unexpected new limits for both ``regular" and singular observables.
	}
\maketitle
\tableofcontents

\section{Introduction}
Digging a hole in the phase space is a powerful tool for analyzing chaotic systems. One classical example dates back to 1970s when Kaplan and Yorke observed an irregular oscillation for the Lorenz system \cite{KaplanYorke79,Lorenz63}. To study it, Yorke and Yorke \cite{YY79} regarded part of the phase space as an escape region (a hole). Pianigiani and Yorke \cite{PY79} formulated such idea mathematically. Since then, the system with holes is known as the open dynamical system. A beautiful survey for open system can be found in \cite{DemersLSY}. 

The Poincar\'e map of Lorenz system exhibits a transition from periodic to chaotic behavior \cite{Manneville1979,PM80}. This intermittent behavior was later intensively studied by many models such as the Liverani--Saussol--Vaienti (LSV) intermittent maps $\{T_{\alpha}, \alpha\in (0,1)\}$. The important feature of these map is the intermittency caused by the neutral fixed point, leading to the loss of spectral gaps and polynomial mixing rates. More details can be found in \cite{LSV99}.

The limit behaviors of closed systems (i.e., systems without holes) attract considerable interest, and many results were established under the SRB measure: for LSV maps, when $0<\alpha<1/2$, Young \cite{Y99} proved Gaussian laws (i.e., central limit theorems) for H\"older observables. When $1/2<\alpha<1$, Melbourne and Zweim\"uller \cite{MZ15} obtained a functional version of stable laws (including Gaussian laws) for H\"older observables. Besides, the Gaussian limit law holds for a large class of exponential mixing systems \cite{Y99,Y,BS}. 

When LSV maps comes with a Markov hole, Demers and Fernandez \cite{DF16} assigned them with a class of initial (non)invariant measure (\ref{h}), and obtained polynomial escape rates (see \cite{S25} for non-Markov holes), and proved that the normalized push-forward measures of the surviving initial distributions converge weakly to the Dirac measure at the neutral fixed point for a subclass of initial measures. 

For some exponential mixing piecewise expanding maps with sufficiently small holes, they exhibit a significantly different phenomenon: the normalized push-forward measures converge to an absolutely continuous conditionally invariant measure (see Bruin, Demers and Melbourne \cite{BDM10}, Lepp\"anen, Nakajima, Nakano \cite{LNN26} and Atnip, Froyland, Gonz\'alez-Tokman, Nakano, Vaienti \cite{AFTNS}). For such systems, conditional limit theorems (i.e., convergence in distribution conditioned on survival orbits up to time $n$) received lots of studies recently, e.g., the work \cite{LNN26,AFTNS} proved that the conditional limit is a Gaussian law. Bassols Cornudella and Castro \cite{BC26} proved that the conditional limit for absorbing Markov chains is a stable law. These conclusions perfectly match the existing limit theory for closed systems. 

However, for the conditional limit theory of open intermittent systems, only some special limit behaviors were numerically investigated by Brevitt and Klages \cite{BK25}. In the present paper, we employ and advance the ideas of the second author's previous work \cite{BunimovichSu}, provide a rigorous conditional limit theory for open LSV maps. That is, we prove
\begin{enumerate}
    \item Theorem \ref{snf_n}: conditional limit for Birkhoff average;
    \item Theorems \ref{snf} and \ref{snf_s}: conditional limit law for a normalized Birkhoff sum;
    \item Theorem \ref{critical_boundedness}: uniform boundedness of conditional expectations.
\end{enumerate}

All these results are new and some are unexpected in the following aspects.
\begin{enumerate}
    \item The initial measure for our limit theory is a broad class of (non)invariant probability absolutely continuous with respect to Lebesgue measure, including SRB/Lebesgue measures, the density function is allowed to oscillate near the neutral fixed point (see more discussions in Remark \ref{measureremark}). This is new compared to \cite{Y99, LNN26} where only  invariant (resp. conditionally invariant) probability was considered for limit theorems for closed (resp. open) systems;
    \item The regularity of our observable is relaxed, e.g., our Theorems \ref{snf_n}-\ref{critical_boundedness} require regularity only at neutral fixed point, and Theorem \ref{snf_s} allows blow-up near neutral fixed point, while H\"older observables were needed in \cite{Y99, MZ15, LNN26} for closed/open systems;
\item Our conditional limit laws (Theorems \ref{snf} and \ref{snf_s}) are not stable laws (see Proposition \ref{nonstable}), let alone the Gaussian laws. This new discovery is unexpected  compared to the existing limit law for closed/open systems \cite{Y99, MZ15, LNN26, AFTNS}. 
\item  The restriction on the hole is relaxed. It can be a disconnected and big hole, not necessarily a Markov hole or a sufficiently small hole required in \cite{DF16,BunimovichSu}.
\item Our Theorem \ref{critical_boundedness} confirms the numeric-based claim for the suppression of Lyapunov stretching in \cite{BK25} (see Remark \ref{physic}).
\item Although only one-dimensional LSV map is studied here, our method can be adapted to a class of multidimensional open intermittent systems \cite{BunimovichSu,highdimlsv} and the same limit distribution still holds (see more discussions in Remark \ref{physic} and Remark \ref{high}).
\end{enumerate}

The paper is organized as follows. Section \ref{notation} introduces notation and conventions. Section \ref{results} states the main results. Section \ref{inducing} includes the inducing scheme of \cite{BunimovichSu} and establishes the return-time and open renewal operator estimates. Section \ref{proofs} proves the main results. Section \ref{appendix} A
ppendix shows that our conditional limit law is not stable, and discusses that the normalization in our limit laws is unique.
    
	\section{Notation and conventions}
	\label{notation}
	\begin{enumerate}
		\item $\mathbb N_0 := \mathbb N \cup \{0\}$.
		
		\item For a function $g$, we write $g(x+):=\lim_{y\downarrow x}g(y)$ whenever this limit exists. 
		
		\item $\mathrm{Leb}$ denotes Lebesgue measure on $[0,1]$. 
		
		\item All observables $f$ and holes $H$ considered in this paper are assumed to be measurable.
		
		\item We write $a_n\sim b_n$ if $a_n/b_n\to1$ as $n\to\infty$.
		
		\item We write $a_n\asymp b_n$ if $cb_n\leq a_n\leq Cb_n$ for some constants $0<c\leq C<\infty$ and all sufficiently large $n$.
		
		\item For $a,b\geq0$, we write $a\lesssim b$ or $a=O(b)$ if $a\leq Cb$ for some constant $C>0$, and $a\lesssim_H b$ or $a=O_H(b)$ if $C$ may depend on $H$.
		
		\item We write $a_n=o(b_n)$ if $a_n/b_n\to0$ as $n\to\infty$.
		
		\item $\mathbb P$ denotes the probability measure on the underlying probability space.

        \item Denote $\inf H$ by $\inf\{x: x\in H \subseteq [0,1]\}$. 
		
		\item We write $X_n \longrightarrow_d X$
		under $\nu_n$, meaning that
		\begin{align*}
			\nu_n \circ X_n^{-1} \longrightarrow_d \mathbb P \circ X^{-1}.
		\end{align*}
		
		\item We write $X_n-Y_n\longrightarrow0$ in $\nu_n$ if for every $\varepsilon>0$,
		\begin{align*}
			\nu_n\bigl(|X_n-Y_n|>\varepsilon\bigr)\longrightarrow0
			\quad\text{as }n\to\infty.
		\end{align*}
	\end{enumerate}

	\section{Main results}
	\label{results}
	Throughout this paper, we consider a family of Liverani--Saussol--Vaienti (LSV) maps $T:[0,1]\to[0,1]$ with $0<\alpha<1$, defined by
	\begin{align*}
		T(x)=
		\begin{cases}
			x(1+2^\alpha x^\alpha), & 0\leq x\leq \frac12,\\[1mm]
			2x-1, & \frac12<x\leq1.
		\end{cases}
	\end{align*}
	Let $H\subsetneq(0,1]$ be a fixed subset (called a hole) bounded away from $0$. Define the first hitting time and the survival set  by
	\begin{align*}
		\tau_H(x):=\inf\{j\geq1:T^jx\in H\},\quad E_n^{H}:=\{x:\tau_H(x)>n\}.
	\end{align*}
	
	Let $\mu$ be an initial probability measure with density $h$. Assume for some $0\leq\beta<1$ and $\varepsilon_0>0$,
	\begin{align}
		\label{h}
		0<\operatorname*{ess\,inf}_{0<x<\varepsilon_0}x^\beta h(x)\leq
		\operatorname*{ess\,sup}_{0<x<1}x^\beta h(x)<\infty.
	\end{align}
	Define the conditional measure on the survival set $E_n^{H}$ by
	\begin{align*}
		\mu_n^{H}(\cdot):=\frac{\mu(\cdot\bigcap E_n^{H})}{\mu(E_n^{H})}.
	\end{align*}
	
	Given an observable $f:[0,1]\to \mathbb{R}$, denote the Birkhoff sums by $S_nf:=\sum_{j=0}^{n-1}f\circ T^j$.
	We now state our first result.
	\begin{theorem}[Conditional limit for Birkhoff average]
		\label{snf_n}
		Let $H\subsetneq(0,1]$ be a fixed hole with
		$\inf H>0$ and $\operatorname{int}H\neq \emptyset$. Assume \eqref{h} with $0\leq\beta<1$. If $f:[0,1]\to\mathbb R$
		is bounded and continuous at $0$, then as $n\to\infty$,
		\begin{align*}
			\frac{S_nf}{n}\longrightarrow f(0)
			\quad \text{in }\mu_n^{H}, \quad \mu_n^{H}\Big(\frac{S_nf}{n}\Big)\longrightarrow f(0).
		\end{align*}
	\end{theorem}

	Next, we study finer statistical properties. Define
	\begin{align}
		\label{jmu}
J_\mu(t):=\mu([0,t])+\sum_{s\geq0}\mu\left(E_s^{H}\cap
		\left\{x:\frac12<T^sx\leq\frac{1+t}{2}\right\}\right).
	\end{align}

	Let $c_\alpha:=(\alpha2^\alpha)^{-1/\alpha}$. For $0\leq\beta<1$ defined in (\ref{h}), let $U$ be a random variable depending on $\alpha, \beta$ such that
	\begin{align*}
		\mathbb P(U>u)=u^{-\frac{1-\beta}{\alpha}}
		\quad\text{for }u\geq1.
	\end{align*}
    
	Given a constant $\gamma \in \mathbb{R}$, define
	\begin{align*}
		G_{\gamma/\alpha}(u)
		:=\int_0^1(u-t)^{-\gamma/\alpha}\,d\mathrm{Leb}(t)
		=\begin{cases}
			\displaystyle
			\frac{u^{1-\gamma/\alpha}-(u-1)^{1-\gamma/\alpha}}{1-\gamma/\alpha},
			&\gamma\neq \alpha,\\[4mm]
			\displaystyle\log\frac{u}{u-1},
			&\gamma=\alpha.
		\end{cases}, \quad u >1.
	\end{align*} Observe that $G_{\gamma/\alpha}$ depends on $\gamma/\alpha$.
	For $\gamma<\alpha$, we extend $G_{\gamma/\alpha}$ continuously to $u=1$
	by setting $G_{\gamma/\alpha}(1):=\alpha/(\alpha-\gamma)$. We now state the conditional limit law for ``regular" observables.
	\begin{theorem}
		\label{snf}
		Let $H\subsetneq(0,1]$ be a fixed hole with
		$\inf H>0$ and $\operatorname{int}H\neq \emptyset$.
		Assume \eqref{h} with $0\leq\beta<1$, let an observable $f:[0,1]\to\mathbb R$ satisfy
		\begin{align*}
			f(x)=f(0)+ax^\gamma+o(x^\gamma)
			\quad\text{as }x\downarrow0
		\end{align*}
		for some $a\neq0$ and $\gamma>0$.
		\begin{enumerate}[label=\textup{(\alph*)}]
			\item Suppose $0<\beta<1$ and for every $\lambda>0$,
			\begin{align}
				\label{mu_rv}
				\frac{\mu([0,\lambda t])}{\mu([0,t])}
				\longrightarrow\lambda^{1-\beta}
				\quad\text{as }t\downarrow0.
			\end{align}
			Then as $n\to\infty$,
			\begin{align*}
				n^{\gamma/\alpha-1}(S_nf-nf(0))
				\longrightarrow_d a c_\alpha^\gamma G_{\gamma/\alpha}(U)
				\quad\text{under }\mu_n^{H}.
			\end{align*}
			
			\item Suppose $\beta=0$, $f$ is bounded, $0<\gamma<\alpha$, and for every $\lambda>0$,
			\begin{align}
				\label{jmu_rv}
				\frac{J_\mu(\lambda t)}{J_\mu(t)}
				\longrightarrow\lambda
				\quad\text{as }t\downarrow0.
			\end{align}
			Then as $n\to\infty$,
			\begin{align*}
				n^{\gamma/\alpha-1}(S_nf-nf(0))
				\longrightarrow_d a c_\alpha^\gamma G_{\gamma/\alpha}(U)
				\quad\text{under }\mu_n^{H}.
			\end{align*}
		\end{enumerate}
	\end{theorem}

When $\beta=0$ and $\gamma\ge \alpha$, we have the following conclusion instead. \begin{theorem}\label{critical_boundedness}
	Let $H\subsetneq(0,1]$ be a fixed hole with $\inf H>0$ and $\operatorname{int}H\neq\emptyset$. Assume \eqref{h} with $\beta=0$. If $f:[0,1]\to\mathbb R$ is bounded and \[f(x)-f(0)=O(x^\gamma) \text{ as } x\downarrow0\] with $\gamma\ge\alpha$, then
	\begin{align*}
		\sup_{n\geq1}\big|\mu_n^H\bigl(S_nf-n f(0)\bigr)\big|\le \sup_{n\geq1}\mu_n^H\bigl(S_n\big|f- f(0)\big|\bigr)<\infty.
	\end{align*}
\end{theorem}

	We next state the limit law for singular observables.
	
	\begin{theorem}
		\label{snf_s}
		Let $H\subsetneq(0,1]$ be a fixed hole with
		$\inf H>0$ and $\operatorname{int}H\neq \emptyset$.
		Assume \eqref{h} with $0\leq\beta<1$, let an observable $f:(0,1]\to\mathbb R$ satisfy
		\begin{align*}
			f(x)=ax^\gamma+o(x^\gamma)
			\quad\text{as }x\downarrow0
		\end{align*}
		for some $a\neq0$ and $\gamma<0$.
		\begin{enumerate}[label=\textup{(\alph*)}]
			\item Suppose $0<\beta<1$ and \eqref{mu_rv} holds.
			Then as $n\to\infty$,
			\begin{align*}
				n^{\gamma/\alpha-1}S_nf
				\longrightarrow_d a c_\alpha^\gamma G_{\gamma/\alpha}(U)
				\quad\text{under }\mu_n^{H}.
			\end{align*}
			
			\item Suppose $\beta=0$, \eqref{jmu_rv} holds,
			and $f$ is bounded on $[\varepsilon,1]$ for any $0<\varepsilon<1$.
			Then as $n\to\infty$,
			\begin{align*}
				n^{\gamma/\alpha-1}S_nf
				\longrightarrow_d a c_\alpha^\gamma G_{\gamma/\alpha}(U)
				\quad\text{under }\mu_n^{H}.
			\end{align*}
		\end{enumerate}
	\end{theorem}

When $\beta=0$, an easy-to-verify condition for (\ref{jmu_rv}) is as follows.
    \begin{corollary}
	\label{snf_bv}
	Let $H$ be a fixed finite union of intervals with $\inf H>0$ and $\operatorname{int}H\neq \emptyset$. Assume \eqref{h} \eqref{mu_rv} with $\beta=0$, and suppose that $h\in\BV([\epsilon,1])$ for any $\epsilon\in (0,1/2)$, i.e., $h$ has bounded variance on $(\epsilon,1]$. 
	\begin{enumerate}[label=\textup{(\alph*)}]
		\item If $f:[0,1]\to\mathbb R$ is bounded, and $f$ satisfies
		\begin{align*}
			f(x)=f(0)+ax^\gamma+o(x^\gamma)
			\quad\text{as }x\downarrow0
		\end{align*}for some $a \neq 0$ and  $0<\gamma<\alpha$, then as $n\to\infty$,
		\begin{align*}
			n^{\gamma/\alpha-1}\bigl(S_nf-nf(0)\bigr)
			\longrightarrow_d a c_\alpha^\gamma G_{\gamma/\alpha}(U)
			\quad \text{under }\mu_n^H.
		\end{align*}
		
		\item If $f:(0,1]\to\mathbb R$ is bounded on $[\varepsilon,1]$ for any $0<\varepsilon<1$, and $f$ satisfies
		\begin{align*}
			f(x)=ax^\gamma+o(x^\gamma)
			\quad\text{as }x\downarrow0
		\end{align*}
		for some $a\neq0$ and  $\gamma<0$, then  as $n\to\infty$, 
		\begin{align*}
			n^{\gamma/\alpha-1}S_nf
			\longrightarrow_d a c_\alpha^\gamma G_{\gamma/\alpha}(U)
			\quad \text{under }\mu_n^H.
		\end{align*}
	\end{enumerate}
\end{corollary}
  \begin{remark}\label{measureremark}
        When $\beta \in [0,1)$, a sufficient condition for \eqref{mu_rv} is
		\begin{align*}
\lim_{t\downarrow0}\frac1t\int_0^t x^{\beta}h(x)\,dx \text{ exists in } (0,\infty).
		\end{align*}
		For example, take $h(x)$ by normalizing functions $x^{-\beta}$ and $
		x^{-\beta}\bigl(2+\sin(1/x)\bigr)$. 
        
        If $h$ is an SRB density function, it follows from \cite{LSV99} that $
		h(x)\sim C_{\alpha} x^{-\alpha}$
		for some $C_{\alpha}>0$.
		Thus, Theorem \ref{snf_n},
		Theorem \ref{snf}\textup{(a)}, and
		Theorem \ref{snf_s}\textup{(a)} hold with $\beta=\alpha\neq0$.

        Our conditional limit law $G_{\gamma/\alpha}(U)$ for open systems depends on $\beta$, the regularity of the density functions of the initial measure. This is new compared to \cite{MZ15, Z} where the limit law does not change for closed systems.
    \end{remark}
	\begin{remark}
The order of the normalization in Theorems \ref{snf} and \ref{snf_s} is not appropriate other than $n^{\gamma/\alpha-1}$. See the discussion in Appendix \ref{rn}.
	\end{remark}

    \begin{remark}
        The observable $f$ in Theorems \ref{snf_n}--\ref{snf_s} can be irregular outside $0$. This is not surprising, because it was already showed in \cite{DemersLSY} that $\mu_n^H \to_d \delta_0$. When $h$ is an SRB density, the second author's work \cite{BunimovichSu} derived a finer quantitative description for this weak convergence.
    \end{remark}
	
	\begin{remark}\label{physic}
	    A direct application of Theorem \ref{critical_boundedness} is that, by taking $f=\log|T'|, \gamma=\alpha, \beta=0$ and $h=1$, we assert suppression of the Lyapunov stretching (see \cite{BK25}) for our LSV maps. The same conclusion holds for Pomeau–Manneville map using the exact same argument, so we rigorously confirm the suppression of the Lyapunov stretching in \cite[table 2 $z\in (1,2)$ and page 13 line -2]{BK25}. 
	\end{remark}
	
	\begin{remark}
		Unlike \cite{BunimovichSu}, our hole $H$ is not necessarily small (same as \cite{S25}). If the hole $H$ contains a neighborhood of the neutral fixed point $0$, the intermittency is killed by the hole, spectral gap exists for the open dynamical system according to \cite{DT17a, DemersLSY}. Hence, following the proof of \cite{AFTNS,LNN26, BC26}, we believe that the conditional limit distribution is stable/Gaussian. 
	\end{remark}
	
	\begin{remark}\label{high}
		Although we only consider one-dimensional intermittent systems, our method is adapted to multidimensional non-Markovian
		nonuniformly expanding systems considered in \cite{BunimovichSu, highdimlsv}, where the intermittency arises from a neutral fixed circle and the initial measure is the SRB probability. The limit distribution for such system is again $G_{\gamma/\alpha}(U)$ up to a constant. For example, if we choose observable $f$ to be bounded and, uniformly in $\theta$,
\[
f(x,\theta)=b+a(\theta)x^\gamma+o(x^\gamma),
\quad x\downarrow0,
\]
where $\gamma>0$, $b\in \mathbb{R}$ and $a(\cdot)$ is H\"older on $\mathbb T$. Then, 
\[
\displaystyle
\frac{S_nf-nb}{n^{1-\gamma/\alpha}}
\longrightarrow_d
c_{\alpha}^\gamma \int_{\mathbb T}a(\theta)\,d\theta \,G_{\gamma/\alpha}(U) \quad \text{under }\mu_n^H.
\] Its proof is similar to Theorem \ref{snf}'s and we will not repeat here.
\end{remark}

	Before giving the proof for our main results, we introduce good holes and some useful tools introduced in \cite{BunimovichSu} in the next section.
	\section{Inducing scheme and renewal operators}
	\label{inducing}
	 In \cite{BunimovichSu}, a general framework for slowing mixing open expanding systems was introduced, inducing method and the open operator renewal theory for small good holes were used when the system preserves an SRB measure. We now summarize them in the case of LSV maps.
	
	Let $T_L(x):=x(1+2^\alpha x^\alpha)$ be the left branch of $T$. Define the preimage sequence by
	\begin{align*}
		x_0:=\frac12,
		\quad
		x_{i+1}:=T_L^{-1}(x_i)
		\quad\text{for }i\geq0.
	\end{align*}
	
	Since $x_i\downarrow0$ and $\inf H>0$, define and fix $
		m:=\min\{i\in\mathbb N_0:x_i<\inf H\}$
	and the inducing base $
		Y:=(x_m,1]$ throughout.
	Particularly, if $\inf H>1/2$, then
	$m=0$ and $Y=(1/2,1]$.
	
	Define the first return time to $Y$ by
	\begin{align*}
		R(x):=\min\{j\geq1:T^jx\in Y\}
	\end{align*}
	and the induced map by
	\begin{align*}
		F:=T^R.
	\end{align*}
	Let $\nu$ be the SRB probability of $F$. According to 
	\cite[Section 6.1.1 and Lemma 22]{BunimovichSu}, the map $F$ is mixing and Gibbs--Markov with big images.
	
	\subsection{Holes and good holes}
	Recall that our hole $H\subsetneq Y$ satisfies 
	$\inf H>0$ and $\operatorname{int}H\neq \emptyset$.
	Define the set of singular points by
	\begin{align*}
		\mathcal S
		:=\bigcup_{j\geq0}T^{-j}\left\{0,\frac12,1\right\}.
	\end{align*}
	An interval $H_g$ is called a
	\emph{good hole} if the following conditions hold:
	
	\begin{itemize}[leftmargin=2em]
		\item The hole $H_g$ has a center
		$z_0\in\operatorname{int}(H_g)\setminus\mathcal S$.
		
		\item The hole $H_g$ satisfies
		$0<\nu(H_g)<\sigma(z_0)$,
		where $\sigma(z_0)>0$ is the threshold given in
		\cite[Theorem 2]{BunimovichSu}.
	\end{itemize}
	
	Obviously, our hole $H$ contains a good hole.
	Indeed, since $\mathcal S$ is countable, we can choose
	$z_0\in(\operatorname{int} H)\setminus\mathcal S$
	and then an interval $H_g$ centered at $z_0$
	such that $H_g\subsetneq \operatorname{int}H$
	and $\nu(H_g)<\sigma(z_0)$. From now on, we fix this $H_g$.
	
	\subsection{Open renewal operators for good holes}
Following
	\cite[Section 6.1.2]{BunimovichSu}, let $\BV$ be the subspace of
	$L^1(Y)$, endowed with the norm
	\begin{align*}
		\|v\|_{\BV}:=|v|_{1}+\operatorname{Var}_Y(v),
	\end{align*}
	where
	\begin{align*}
		\operatorname{Var}_Y(v)
		:=\inf_{\widetilde v=v\;m_Y\text{-a.e.}}
		\sup_{\substack{N\ge1\\ x_m<t_0<t_1<\cdots<t_N\le1}}
		\sum_{i=1}^{N}\bigl|\widetilde v(t_i)-\widetilde v(t_{i-1})\bigr|.
	\end{align*}
		Let $\mathcal L$ be the transfer operator of $F$ with respect to Lebesgue measure. Define
	\begin{align*}
		\mathcal R_k v:=\mathcal L\bigl(\1_{\{R=k\}}v\bigr)
		\quad\text{for all }k\geq1.
	\end{align*}
	Then \cite[Lemma~21]{BunimovichSu} gives
	\begin{align*}
		\|\mathcal R_k\|_{\BV\to\BV}
		\lesssim_H k^{-1/\alpha-1}.
	\end{align*}
	
	We now introduce open transfer operators for the good hole $H_g$. Define
	\begin{align*}
		\mathring{\mathcal L}(\cdot):=\1_{Y\setminus H_g}\mathcal L(\cdot),\quad 
		\mathring{\mathcal R}_k(\cdot):=\mathring{\mathcal L}\bigl(\1_{\{R=k\}}\cdot \bigr),\quad 
		\mathring{\mathcal R}(z)
		:=
		\sum_{k\geq1}z^k\mathring{\mathcal R}_k.
	\end{align*}

    It follows from \cite{BunimovichSu} that the spectral radius of $\mathring{\mathcal R}(z)$ is no more than $|z|$. Hence for any $|z|<1$, we can write
		\begin{align*}
			[I-\mathring{\mathcal R}(z)]^{-1}
			:=\sum_{n\geq0}\mathring{\mathcal T}_nz^n.
		\end{align*} where $\mathring{\mathcal T}_0=Id$ and
        \begin{align*}
			\mathring{\mathcal T}_n
			=\sum_{\ell=1}^n
			\sum_{\substack{k_1,\ldots,k_\ell\geq1\\
					k_1+\cdots+k_\ell=n}}
			\mathring{\mathcal R}_{k_\ell}\circ\cdots\circ
			\mathring{\mathcal R}_{k_1}.
		\end{align*} They have the following estimate.
	\begin{proposition}
		\label{tn}There is a constant $C>0$ such that for any $n\ge 1$, \[
			\|\mathring{\mathcal T}_n\|_{\BV\to\BV}
			\le C(n+1)^{-1-1/\alpha}.
		\] Besides, $\mathring{\mathcal T}_n$ can be extended to $L^1$, and for any bounded function $v$ (might not be of bounded variance),
		\begin{align*}
			|\mathring{\mathcal T}_nv|_\infty
			\le C(n+1)^{-1-1/\alpha}|v|_\infty.
		\end{align*}
	\end{proposition}
	
	\begin{proof}
		The estimate
		\begin{align*}
			\|\mathring{\mathcal T}_n\|_{\BV\to\BV}
			\le C(n+1)^{-1-1/\alpha}
		\end{align*}
		follows from \cite[Lemma~12]{BunimovichSu}. For any $n\geq1$,
		\begin{align*}
			\mathring{\mathcal T}_n
			=\sum_{\ell=1}^n
			\sum_{\substack{k_1,\ldots,k_\ell\geq1\\
					k_1+\cdots+k_\ell=n}}
			\mathring{\mathcal R}_{k_\ell}\circ\cdots\circ
			\mathring{\mathcal R}_{k_1}.
		\end{align*}Since each $\mathring{\mathcal R}_k$ maps $L^1$ to $L^1$, so $\mathring{\mathcal T}_n$ can be extended to $L^1$. 
		Since each $\mathring{\mathcal R}_k$ is a positive operator, $\mathring{\mathcal T}_n$ is positive too. 
		Moreover,
		$0\leq |v|\leq |v|_\infty\1$ implies
		\begin{align*}
			|\mathring{\mathcal T}_nv|
			\leq\mathring{\mathcal T}_n(|v|)
			\leq |v|_\infty\mathring{\mathcal T}_n\1.
		\end{align*}
		Using
		$\BV(Y)\hookrightarrow L^\infty(Y)$, we have
		\begin{align*}
			|\mathring{\mathcal T}_nv|_\infty
			\leq |v|_\infty |\mathring{\mathcal T}_n\1|_\infty
			\lesssim |v|_\infty\|\mathring{\mathcal T}_n\1\|_{\BV}
			\lesssim_{H_g}(n+1)^{-1-1/\alpha} |v|_\infty.
		\end{align*}
	\end{proof}

    \subsection{Tail asymptotics}
	We prove finer asymptotic estimates for $x_n$, which will be used in the proofs of our main theorems.
	\begin{proposition}
		\label{xi_xi-1}
		For any $i>0$, any $x_i<x\leq x_{i-1}$,
		\begin{align*}
			x_i
			=c_\alpha i^{-1/\alpha}
			\left(1+O\left(\frac{\log i}{i}\right)\right),\quad 
			x=c_\alpha i^{-1/\alpha}\left(1+O\left(\frac{\log i}{i}\right)\right).
		\end{align*}
        where the constants in $o(\cdot),O(\cdot)$ do not depend on $i,x$.
	\end{proposition}
	
	\begin{proof}
		Let $u_i:=x_i^{-\alpha}$. By
		\begin{align*}
			x_{i+1}^{-\alpha}-x_i^{-\alpha}
			=x_{i+1}^{-\alpha}\left[1-\left(1+2^\alpha x_{i+1}^\alpha\right)^{-\alpha}\right],
		\end{align*}
		we have
		\begin{align*}
			u_{i+1}-u_i
			=u_{i+1}\left[1-\left(1+\frac{2^\alpha}{u_{i+1}}\right)^{-\alpha}\right]
			=u_{i+1}\left(\frac{\alpha 2^\alpha}{u_{i+1}}+O\bigl(u_{i+1}^{-2}\bigr)\right)
			=\alpha 2^\alpha+O\left({1}/{u_{i+1}}\right).
		\end{align*}
		Hence $u_i\asymp i$. As $i\to\infty$, ${1}/{u_{i+1}}=O(1/i)$ and
		\begin{align*}
			u_{i+1}-u_i
			=\alpha2^\alpha+O\left(1/i\right),
		\end{align*}
		so
		\begin{align*}
			u_i
			=
			u_0+\sum_{j=0}^{i-1}(u_{j+1}-u_j)=\alpha2^\alpha i+O(\log i).
		\end{align*}
		It follows that
		\begin{align*}
			x_i
			=(\alpha2^\alpha)^{-1/\alpha}
			i^{-1/\alpha}
			\left(1+O\left(\frac{\log i}{i}\right)\right).
		\end{align*}
		
		Using the above estimate, we have
		\begin{align*}
			x_i&=c_\alpha i^{-1/\alpha}\left(1+O\left(\frac{\log i}{i}\right)\right),\\
			x_{i-1}&=c_\alpha (i-1)^{-1/\alpha}\left(1+O\left(\frac{\log i}{i}\right)\right).
		\end{align*}
		Hence for $x_i<x\leq x_{i-1}$,
		\begin{align*}
			x
			=
			c_\alpha i^{-1/\alpha}\left(1+O\left(\frac{\log i}{i}\right)\right).
		\end{align*}
	\end{proof}

	\begin{proposition}
		\label{tl_j}
		For any $x\in (x_m,T_L(x_m)]$,
		\begin{align*}
			T_L^{-j}(x)\asymp j^{-1/\alpha},\quad 
			(T_L^{-j})'(x)\asymp j^{-1-1/\alpha}
		\end{align*}where the constants in $\asymp$ do not depend on $j>0,x$.
	\end{proposition}
	
	\begin{proof}
		For $x\in(x_m,T_L(x_m)]$, we have
		\begin{align*}
			x_{m+j}<T_L^{-j}(x)\leq x_{m+j-1}.
		\end{align*}
		Then
		Proposition~\ref{xi_xi-1} gives
		\begin{align*}
			T_L^{-j}(x)=c_\alpha j^{-1/\alpha}
			\left(1+O({\log j}/{j})\right).
		\end{align*}
		
		We now estimate the derivative $(T_L^{-j})'(x)$.
		By $
		T_L'(x)=1+(1+\alpha)2^{\alpha}x^{\alpha},$
		we have
		\begin{align*}
			(T_L^{-j})'(x)
			&=\prod_{k=1}^{j}\frac{1}{T_L'\bigl(T_L^{-k}(x)\bigr)}
			=\prod_{k=1}^{j}
			\frac{1}{1+(1+\alpha)2^{\alpha}\bigl(T_L^{-k}(x)\bigr)^{\alpha}},
		\end{align*}
		Then
		\begin{align*}
			\log\bigl((T_L^{-j})'(x)\bigr)
			&=-\sum_{k=1}^{j}
			\log\left(1+(1+\alpha)2^{\alpha}
			\bigl(T_L^{-k}(x)\bigr)^{\alpha}\right).
		\end{align*}
		Besides, we have
		\begin{align*}
			2^{\alpha}\bigl(T_L^{-k}(x)\bigr)^{\alpha}
			=2^{\alpha}c_{\alpha}^{\alpha}k^{-1}
			\left(1+O\left(\frac{\log k}{k}\right)\right)^{\alpha}
			=\frac{1}{\alpha k}
			+O\left(\frac{\log(k+1)}{k^{2}}\right),
		\end{align*}
		where we use $
		2^{\alpha}c_{\alpha}^{\alpha}=1/\alpha.$
		Using \(\log(1+t)=t+O(t^2)\) as \(t\to 0\), we get
		\[
		\log\left(1+(1+\alpha)2^{\alpha}
		\bigl(T_L^{-k}(x)\bigr)^{\alpha}\right)
		=
		\left(1+\frac{1}{\alpha}\right)\frac{1}{k}
		+O\left(\frac{\log(k+1)}{k^{2}}\right).
		\]
		Hence,
		\begin{align*}
			\log\bigl((T_L^{-j})'(x)\bigr)
			=-\sum_{k=1}^{j}
			\log\left(1+(1+\alpha)2^{\alpha}
			\bigl(T_L^{-k}(x)\bigr)^{\alpha}\right)
			=-\left(1+\frac{1}{\alpha}\right)\log j+O(1).
		\end{align*}
		Remove the $\log$ we conclude the proof.
	\end{proof}

	\section{Proofs of the main results}
	\label{proofs}	Inspired by the ideas of the proof in \cite{BunimovichSu} where it derived an asymptotic expansion for the measure of the complicated set $E_n^H$ based on a decomposition of $E_n^H$, we employ the same strategy: find a suitable decomposition of $E_n^H$, then apply it to prove Theorems \ref{snf_n},
	\ref{snf}, \ref{critical_boundedness}, and \ref{snf_s},
	and Corollary \ref{snf_bv}.
    
	\subsection{Decomposition of \texorpdfstring{$E_n^{H}$}{E-n-H}}
	For each $n\geq1$, we divide $E_n^H$ into three disjoint classes:
	\begin{enumerate}
		\item Points starting in $[0,x_{m+n}]$.
		\item Points ending in $[0,x_m]$ at time $n$, whose last entry into this interval from $Y$ occurs at some time $s\in\{1,\ldots,n\}$. Their position at time $s$ lies in $[0,x_{m+n-s}]$.
		\item Points ending in $Y$ at time $n$.
	\end{enumerate}
	
	To record the last entry time, for $n\geq0$ let $\sigma_n(x)$ be the largest $s\in\{1,\ldots,n\}$ such that $T^{s-1}x\in Y$ and $T^s x\in[0,x_m]$. Set $\sigma_n(x)=0$ if no such time exists. On \(\{T^nx\in[0,x_m]\}\), the cases \(\sigma_n=0\) and \(\sigma_n\ge1\) correspond to the first and second classes, respectively. We now estimate the measure of these three classes of sets.
	
	\begin{proposition}
		\label{en}
		Under \eqref{h}, for any $n\geq1$, we have
		\begin{align*}
			\mu([0,x_{m+n}])&\asymp n^{-\frac{1-\beta}{\alpha}},\\
			\mu(E_n^{H}\cap\{T^nx\in[0,x_m],\ \sigma_n\geq1\})
			&\lesssim n^{-\frac1\alpha}+n^{-\frac{1-\beta}{\alpha}-1},
			\\
			\mu(E_n^{H}\cap\{T^nx\in Y\})
			&\lesssim n^{-\frac{1-\beta}{\alpha}-1},
		\end{align*}where the constants in $\lesssim,\asymp$ do not depend on $n$.
	\end{proposition}
	
	\begin{remark}\label{1and2overlap}
		Obviously, for $0<\beta<1$, $\mu ([0,x_{m+n}])\asymp \mu(E_n^H)$. For $\beta=0$,  the second class may contribute the first class at the same order $n^{-1/\alpha}$, but both dominate the third class.
	\end{remark}
	
	We define the surviving density at time $n$ by
	\begin{align*}
		h_n^{H}:=\frac{d\,T_*^n(\mu|_{E_n^{H}})}{d\mathrm{Leb}}.
	\end{align*}
	Then $h_0^{H}=h$ and
	\begin{align*}
		\int_0^1h_n^{H}(x)\,d\mathrm{Leb}(x)=\mu(E_n^{H}).
	\end{align*}
	For $s\geq1$, define 
	\begin{align*}
		e_s^{H}(x):=\frac12 h_{s-1}^{H}\left(\frac{1+x}{2}\right)
		 \text{ for any } x\in[0,x_m]
	\end{align*}
	and let $e_0^{H}(x):=h(x)$ for any $x\in [0,x_m]$. To prove Proposition \ref{en}, we need the following two propositions.
	
	\begin{proposition}
		\label{hn_es}
		Under \eqref{h} with $0\leq\beta<1$, for all $n,s\geq1$,
		\begin{align*}
			\|h_n^{H}\|_{L^\infty(Y)}
			\lesssim n^{-\frac{1-\beta}{\alpha}-1},\qquad 
			\|e_s^{H}\|_{L^\infty([0,x_m])}
			\lesssim s^{-\frac{1-\beta}{\alpha}-1}
		\end{align*}where the constants $\lesssim,\asymp$ do not depend on $n,s$.
	\end{proposition}
	
	\begin{proof}
		First, we consider $h_n^{H}|_Y$. For the good hole, write
		$E_n^{H_g}:=\bigcap_{j=1}^nT^{-j}(H_g^c)$
		and
		\begin{align*}
			h_n^{H_g} := \frac{d\,T_*^n\big(\mu|_{E_n^{H_g}}\big)}{d\mathrm{Leb}}.
		\end{align*}
		Using $H_g\subsetneq H$, we split surviving orbits into the ones starting from $Y$ then ending in $Y$, and the ones first entering $Y$ at time $j\geq1$ then ending in $Y$. This gives
		\begin{align*}
			h_n^{H}|_Y
			&\leq h_n^{H_g}|_Y
			\leq\mathring{\mathcal T}_n(h|_Y)
			+\sum_{j=1}^n\mathring{\mathcal T}_{n-j}
			\left[\1_{(x_m,T_L(x_m)]}(h\circ T_L^{-j})(T_L^{-j})'\right].
		\end{align*}
		For each term in $\sum_j^n$, 
		Proposition \ref{tl_j} gives
		\begin{align*}
			h(T_L^{-j}x)(T_L^{-j})'(x)
			\leq\left(\operatorname*{ess\,sup}_{0<x<1}x^\beta h(x)\right)
			(T_L^{-j}x)^{-\beta}(T_L^{-j})'(x)
			\lesssim_{H,\mu}j^{\beta/\alpha}j^{-1-1/\alpha}
			=j^{-\frac{1-\beta}{\alpha}-1}.
		\end{align*}
		For orbits starting in $Y$, we have
		\begin{align*}
			\|h\|_{L^\infty(Y)}
			\leq x_m^{-\beta}\operatorname*{ess\,sup}_{0<x<1}x^\beta h(x)
			<\infty.
		\end{align*}
		Combining these estimates with Proposition \ref{tn}, we obtain
		\begin{align*}
			\|h_n^{H}\|_{L^\infty(Y)}
			&\lesssim_{H,\mu}|\mathring{\mathcal T}_n\1|_\infty
			+\sum_{j=1}^n j^{-\frac{1-\beta}{\alpha}-1}|\mathring{\mathcal T}_{n-j}\1|_\infty\\
			&\lesssim_{H,\mu}n^{-\frac{1}{\alpha}-1}
			+\sum_{j=1}^n j^{-\frac{1-\beta}{\alpha}-1}(n-j+1)^{-\frac{1}{\alpha}-1}\\
            &\lesssim_{H,\mu}n^{-\frac{1}{\alpha}-1}
			+n^{-\frac{1-\beta}{\alpha}-1}
			\lesssim_{H,\mu}n^{-\frac{1-\beta}{\alpha}-1}.
		\end{align*}
		Hence $\|h_n^{H}\|_{L^\infty(Y)}
			\lesssim_{H,\mu}n^{-\frac{1-\beta}{\alpha}-1}$. The estimate for $\|e_s^{H}\|_{L^\infty([0,x_m])}$ follows from its definition.
	\end{proof}
	
	\begin{proposition}
		\label{sigma_un}
		Under \eqref{h} with $0\leq\beta<1$, for all $n\geq1$ and $0\leq K<n$,
		\begin{align*}
			\mu_n^{H}(T^nx\in Y)\lesssim n^{-1},\quad
			\mu_n^{H}(\sigma_n>K)
			\lesssim (K+1)^{-\frac{1-\beta}{\alpha}}+n^{-1}.
		\end{align*}where the constants in $\lesssim$ do not depend on $K,n$.
	\end{proposition}
	
	\begin{proof}
		Since $[0,x_{m+n}]\subseteq E_n^{H}$, then
		\begin{align*}
			\mu(E_n^{H})
			&\geq\mu([0,x_{m+n}])
			=\int_0^{x_{m+n}}h(x)\,d\mathrm{Leb}(x)
			\asymp {x_{m+n}^{1-\beta}}
			\asymp n^{-\frac{1-\beta}{\alpha}}.
		\end{align*}
        
		Using Proposition \ref{hn_es}, we have
		\begin{align*}
			\mu_n^{H}(T^nx\in Y)=\frac{\int_Yh_n^{H}(x)\,d\mathrm{Leb}(x)}{\mu(E_n^{H})}
			\lesssim\frac{\|h_n^{H}\|_{L^\infty(Y)}}{\mu(E_n^{H})}\lesssim \frac{n^{-\frac{1-\beta}{\alpha}-1}}{n^{-\frac{1-\beta}{\alpha}}}
			=n^{-1}.
		\end{align*}

		We now estimate the conditional tail of $\sigma_n$. Since
		\begin{align*}
			\mu_n^{H}(\sigma_n>K)
			&=\mu_n^{H}(\sigma_n>K,\ T^nx\in[0,x_m])
			+\mu_n^{H}(\sigma_n>K,\ T^nx\in Y)\\
			&\leq
			\frac{\mu(E_n^{H}\cap\{T^nx\in[0,x_m],\ \sigma_n>K\})}{\mu(E_n^{H})}
			+\mu_n^{H}(T^nx\in Y)\\
			&\lesssim \frac{\displaystyle\sum_{s=K+1}^n
				\mu(E_n^{H}\cap\{T^nx\in[0,x_m],\ \sigma_n=s\})}{\mu(E_n^{H})}
			+n^{-1}.
		\end{align*}
		
        Now we estimate
		\begin{align*}
			\mu(E_n^{H}\cap\{T^nx\in[0,x_m],\ \sigma_n=s\})
			\text{ for all }0\leq s\leq n.
		\end{align*}
		
		First we claim that
		\begin{align}
			\label{es}
			\mu(E_n^{H}\cap\{T^nx\in[0,x_m],\ \sigma_n=s\})
			=\int_0^{x_{m+n-s}}e_s^{H}(x)\,d\mathrm{Leb}(x).
		\end{align}
		 Indeed, if $s=0$, we have
		\begin{align*}
			\mu(E_n^{H}\cap\{T^nx\in[0,x_m],\ \sigma_n=0\})
			&=\mu([0,x_{m+n}])\\
			&=\int_0^{x_{m+n}}h(x)\,d\mathrm{Leb}(x)
			=\int_0^{x_{m+n}}e_0^{H}(x)\,d\mathrm{Leb}(x).
		\end{align*}
		If $1\leq s\leq n$,
		\begin{align*}
			\mu(E_n^{H}\cap\{T^nx\in[0,x_m],\ \sigma_n=s\})
			&=\mu\left(E_{s-1}^{H}\cap
			\left\{\frac12<T^{s-1}x\leq\frac{1+x_{m+n-s}}2\right\}\right)\\
			&=\int_{1/2}^{(1+x_{m+n-s})/2}h_{s-1}^{H}(x)\,d\mathrm{Leb}(x)\\
			&=\int_0^{x_{m+n-s}}e_s^{H}(x)\,d\mathrm{Leb}(x).
		\end{align*}
		Combining Proposition \ref{hn_es} and Proposition \ref{xi_xi-1},
		we have
		\begin{align}
			\sum_{s=K+1}^n\int_0^{x_{m+n-s}}
			e_s^{H}(x)\,d\mathrm{Leb}(x)
			&\lesssim_{H,\mu}\sum_{s=K+1}^n
			(n-s+1)^{-\frac{1}{\alpha}}s^{-\frac{1-\beta}{\alpha}-1}\notag\\
			&=\left[\sum_{K<s\leq n/2}+\sum_{\max\{K,n/2\}<s\leq n}\right]
			(n-s+1)^{-\frac{1}{\alpha}}s^{-\frac{1-\beta}{\alpha}-1}\notag\\
			&\lesssim n^{-\frac{1}{\alpha}}\sum_{s=K+1}^\infty s^{-\frac{1-\beta}{\alpha}-1}
			+n^{-\frac{1-\beta}{\alpha}-1}\sum_{j\geq1}j^{-\frac{1}{\alpha}}\notag\\
			&\lesssim 
			n^{-\frac{1}{\alpha}}(K+1)^{-\frac{1-\beta}{\alpha}}+n^{-\frac{1-\beta}{\alpha}-1}.
			\label{es_k}
		\end{align}
		
		Thus, for all $0\leq K<n$,
		\eqref{es} and
		\eqref{es_k} give
		\begin{align*}
			\mu_n^{H}(\sigma_n>K)
            \lesssim  (K+1)^{-\frac{1-\beta}{\alpha}}+n^{-1}.
		\end{align*}
	\end{proof}
	
	\begin{lemma}\label{sigman/n}
		As $n\to\infty$,
		\begin{align*}
			\frac{\sigma_n}{n}\longrightarrow0
			\quad \text{in }\mu_n^{H},
		\end{align*}which confirms the simulation observation in \cite[page 9, line 16]{BK25}.
	\end{lemma}
	\begin{proof}
	For every $0<\varepsilon<1$,
	by Proposition \ref{sigma_un} with
	$K=\lfloor\varepsilon n\rfloor$,
	\begin{align*}
		\mu_n^{H}\left(\frac{\sigma_n}{n}>\varepsilon\right)
		=\mu_n^{H}(\sigma_n>\lfloor\varepsilon n\rfloor)\lesssim_{H,\mu}(\lfloor\varepsilon n\rfloor+1)^{-\frac{1-\beta}{\alpha}}+n^{-1}
		\longrightarrow0.
	\end{align*}
    \end{proof}
	
	Proposition \ref{en} follows from Propositions \ref{hn_es}
	and \ref{sigma_un}.
	\begin{proof}[Proof of Proposition~\ref{en}]
		For the first class, we have
		\begin{align*}
			\mu([0,x_{m+n}])
			\asymp x_{m+n}^{1-\beta}
			\asymp n^{-\frac{1-\beta}{\alpha}}.
		\end{align*}
		
		For the second class, we partition according to $\sigma_n=s$. Taking $K=0$ in \eqref{es_k} gives
		\begin{align*}
			\mu(E_n^{H}\cap\{T^nx\in[0,x_m],\ \sigma_n\geq1\})
			=\sum_{s=1}^n\int_0^{x_{m+n-s}}e_s^{H}(x)\,d\mathrm{Leb}(x)\lesssim_{H,\mu}n^{-\frac1\alpha}+n^{-\frac{1-\beta}{\alpha}-1}.
		\end{align*}
		
		For the third class, Proposition~\ref{hn_es} gives
		\begin{align*}
			\mu(E_n^{H}\cap\{T^nx\in Y\})
			=\int_Yh_n^{H}(x)\,d\mathrm{Leb}(x)
			\leq\|h_n^{H}\|_{L^\infty(Y)}
			\lesssim_{H,\mu}n^{-\frac{1-\beta}{\alpha}-1}.
		\end{align*}
	\end{proof}
	
	\subsection{Proof of Theorem \ref{snf_n} }
	In this subsection, we prove the conditional limit for Birkhoff average using Proposition \ref{sigma_un} and the continuity of $f$ at $0$.
	\begin{proof}[Proof of Theorem \ref{snf_n}]
		Fix $\eta>0$. Choose
		$0<\varepsilon<x_m$ such that
		\begin{align*}
			\omega_f(\varepsilon)
			:=\sup_{0\leq x\leq\varepsilon}|f(x)-f(0)|<\frac\eta2.
		\end{align*}
		Define
		\begin{align*}
			A_{n,\eta}:=\left\{x:
			\left|\frac{S_nf(x)}n-f(0)\right|>\eta\right\}.
		\end{align*}
		Then by Proposition \ref{sigma_un}, 
		\begin{align}
			\mu_n^{H}(A_{n,\eta})
			&=\mu_n^{H}(A_{n,\eta}\cap\{T^nx\in Y\})+\mu_n^{H}(A_{n,\eta}\cap\{T^nx\in[0,x_m]\})\notag\\
			&\lesssim n^{-1}+\mu_n^{H}(A_{n,\eta}\cap\{T^nx\in[0,x_m]\}).
			\label{muAneta}
		\end{align}
			
		First, we claim that on $\{T^nx\in[0,x_m]\}$,
		\begin{align}
			\label{snf_w}
			\left|\frac{S_nf(x)}n-f(0)\right|
			\leq\omega_f(\varepsilon)
			+\frac{2\|f\|_\infty(\sigma_n+\frac{x_m}{2^\alpha\varepsilon^{1+\alpha}})}n.
		\end{align}
		Indeed,
		\begin{align*}
			\left|\frac{S_nf(x)}n-f(0)\right|
			&=\frac1n\left|\sum_{j=0}^{n-1}
			\bigl(f(T^jx)-f(0)\bigr)\right|\\
			&\leq\frac1n\sum_{j=0}^{\sigma_n-1}|f(T^jx)-f(0)|
			+\frac1n\sum_{j=\sigma_n}^{n-1}|f(T^jx)-f(0)|\\
			&=\frac1n\sum_{j=0}^{\sigma_n-1}|f(T^jx)-f(0)|
			+\frac1n \Bigl(
			\sum_{\substack{\sigma_n\leq j<n\\T^jx<\varepsilon}}
			+
			\sum_{\substack{\sigma_n\leq j<n\\T^jx\geq\varepsilon}}
			\Bigr)
			|f(T^jx)-f(0)|.
		\end{align*}
		The boundedness of $f$ gives
		\begin{align*}
			\frac1n\sum_{j=0}^{\sigma_n-1}|f(T^jx)-f(0)|
			\leq\frac{2\|f\|_\infty\sigma_n}{n},
		\end{align*}
		and the definition of $\omega_f(\varepsilon)$ gives
		\begin{align*}
			\frac1n\sum_{\substack{\sigma_n\leq j<n\\T^jx<\varepsilon}}
			|f(T^jx)-f(0)|
			\leq\frac{n-\sigma_n}{n}\omega_f(\varepsilon)
			\leq\omega_f(\varepsilon).
		\end{align*}
		It remains to estimate the sum over $T^jx\geq\varepsilon$. Let
		\begin{align*}
			N_{n,\varepsilon}(x)
			:=\sum_{j=\sigma_n}^{n-1}\1_{[\varepsilon,x_m]}(T^jx).
		\end{align*}
		Since $T^jx\in[0,x_m]$ for $\sigma_n\leq j\leq n$, we have
		\begin{align*}
			\frac1n\sum_{\substack{\sigma_n\leq j<n\\T^jx\geq\varepsilon}}
			|f(T^jx)-f(0)|
			\leq\frac{2\|f\|_\infty N_{n,\varepsilon}(x)}n.
		\end{align*}
		Combining these estimates, we obtain
		\begin{align*}
			\left|\frac{S_nf(x)}n-f(0)\right|
			\leq\omega_f(\varepsilon)
			+\frac{2\|f\|_\infty}{n}
			\bigl(\sigma_n+N_{n,\varepsilon}(x)\bigr).
		\end{align*}
		We now bound $N_{n,\varepsilon}$ to prove \eqref{snf_w}. For $\sigma_n\leq j<n$, we have
		\begin{align*}
			T^{j+1}x-T^jx
			=T_L(T^jx)-T^jx
			=2^\alpha(T^jx)^{1+\alpha}.
		\end{align*}
		Hence
		\begin{align*}
			N_{n,\varepsilon}(x)
			&=\sum_{j=\sigma_n}^{n-1}\1_{[\varepsilon,x_m]}(T^jx)\leq\frac1{2^\alpha\varepsilon^{1+\alpha}}
			\sum_{j=\sigma_n}^{n-1}
			2^\alpha(T^jx)^{1+\alpha}\1_{[\varepsilon,x_m]}(T^jx)\\
			&=\frac1{2^\alpha\varepsilon^{1+\alpha}}
			\sum_{j=\sigma_n}^{n-1}
			\bigl(T^{j+1}x-T^jx\bigr)\1_{[\varepsilon,x_m]}(T^jx)\leq\frac1{2^\alpha\varepsilon^{1+\alpha}}
			\sum_{j=\sigma_n}^{n-1}\bigl(T^{j+1}x-T^jx\bigr)\\
			&
			\leq\frac{x_m}{2^\alpha\varepsilon^{1+\alpha}},
		\end{align*}
		so \eqref{snf_w} holds. Thus
		\begin{align*}
			\mu_n^{H}(A_{n,\eta}\cap\{T^nx\in[0,x_m]\}) 
			\leq
			\mu_n^{H}\left(
			\omega_f(\varepsilon)
			+\frac{2\|f\|_\infty
				\left(\sigma_n+\frac{x_m}{2^\alpha\varepsilon^{1+\alpha}}\right)}n>\eta
			\right).
		\end{align*}
		Combining this with \eqref{muAneta} and Lemma \ref{sigman/n}, we obtain
		\begin{align*}
			\mu_n^{H}(A_{n,\eta})
			&\lesssim n^{-1}
			+\mu_n^{H}\left(
			\omega_f(\varepsilon)
			+\frac{2\|f\|_\infty
				\left(\sigma_n+\frac{x_m}{2^\alpha\varepsilon^{1+\alpha}}\right)}n>\eta
			\right)\\
            &\lesssim n^{-1}
			+\mu_n^{H}\left(
			\frac{2\|f\|_\infty
				\left(\sigma_n+\frac{x_m}{2^\alpha\varepsilon^{1+\alpha}}\right)}n>\eta/2
			\right)\to 0, 
		\end{align*}which concludes the proof of conditional convergence in probability.
        
        By the boundedness of $f$, let $n\to \infty$ then $\eta \to 0$, we have
        \begin{align*}
            \Big|\mu_n^{H}\Big(\frac{S_nf}{n}\Big)-f(0)\Big|=\Big|\int_{A_{n,\eta}} \frac{S_nf}{n}-f(0)d\mu_{n}^H\Big|+O(\eta) \to 0 
        \end{align*}which concludes the proof of conditional convergence in expectation.
	\end{proof}
	
	\subsection{Proof of Theorem \ref{snf}}
	For $n\geq1$ and $x\in(0,1]$, define
	\begin{align*}
		U_n(x):=\frac1n\inf\{j\geq n:T^jx\in Y\}.
	\end{align*}

    \subsubsection{Scheme of proof}
\begin{enumerate}
    \item \label{1} When $0<\beta<1$, first apply Proposition \ref{en} to show $\mu_n^{H}(\cdot)
			$ can be approximated by $ \mu_n^{H}(\,\cdot\mid[0,x_{m+n}])$, where the dynamics of the survival orbits on $[0,x_{m+n}]$ is much simpler.
            \item  \label{2} Show that $U_n\longrightarrow_d U$ under $\mu_n^H(\cdot|[0,x_{m+n}])$.
            \item \label{3} Show that, in a space with a large probability, $n^{\gamma/\alpha-1}\bigl(S_nf-nf(0)\bigr)$ can be approximated by $a c_\alpha^\gamma  G_n(U_n)$ for some function $G_n$.
            \item \label{4} In a space with a large probability, prove $G_n(U_n) \asymp G_{\gamma/\alpha}(U_n)$ pointwise.
            \item \label{5} Prove the limit distribution $ ac_{\alpha}^{\gamma}G_{\gamma/\alpha}(U)$. 
    \item  The case $\beta=0$ follows similar strategy.
\end{enumerate} 

In what follows, these statements will be clear.
	
	\subsubsection{Proof of Theorem~\ref{snf}\textup{(a)}}
	\subsubsection*{Step \ref{1}: approximation by $[0,x_{m+n}]$}
	\begin{lemma}
		\label{en_tv}
		As $n\to\infty$,
		\begin{align*}
			\left\|\mu_n^{H}(\cdot)
			-\mu_n^{H}(\,\cdot\mid[0,x_{m+n}])\right\|_{\TV}
			\lesssim n^{-\beta/\alpha}+n^{-1}\longrightarrow0
		\end{align*}where the constant in $\lesssim$ does not depend on $n$.
	\end{lemma}
	
	\begin{proof}
		Since $[0,x_{m+n}] \subsetneq E_n^{H}$,	by Proposition \ref{en},
		\begin{align*}
			\mu(E_n^{H}\setminus[0,x_{m+n}])
			&\lesssim n^{-\frac1\alpha}+n^{-\frac{1-\beta}{\alpha}-1},\\
			\mu(E_n^{H})&\geq\mu([0,x_{m+n}])
			\asymp n^{-\frac{1-\beta}{\alpha}}.
		\end{align*}
	Hence
		\begin{align*}
			\left\|\mu_n^{H}(\cdot)
			-\mu_n^{H}(\,\cdot\mid[0,x_{m+n}])\right\|_{\TV}&\lesssim  2 \frac{\mu(E_n^{H}\setminus[0,x_{m+n}])}{\mu(E_n^{H})}\\
            &\lesssim \frac{n^{-\frac{1}{\alpha}}+n^{-\frac{1-\beta}{\alpha}-1}}{n^{-\frac{1-\beta}{\alpha}}}
			\lesssim n^{-\beta/\alpha}+n^{-1}\longrightarrow0.
		\end{align*}
	\end{proof}
    \subsubsection*{Step \ref{2}: convergence of $U_n$}	
    \begin{lemma}
		\label{un}
		As $n\to\infty$,
		\begin{align*}
			U_n\longrightarrow_d U
			\quad \text{under }\mu_n^{H}(\,\cdot\mid[0,x_{m+n}]).
		\end{align*}
	\end{lemma}
	
	\begin{proof}
		For $u\geq1$, we have
		\begin{align*}
			\mu_n^{H}(U_n>u\mid[0,x_{m+n}])
			=\frac{\mu_n^{H}([0,x_{m+\lfloor un\rfloor}])}
			{\mu_n^{H}([0,x_{m+n}])}
			=\frac{\mu([0,x_{m+\lfloor un\rfloor}])}{\mu([0,x_{m+n}])}.
		\end{align*}
		By Proposition~\ref{xi_xi-1}, as $n\to\infty$,
		\begin{align*}
			\frac{x_{m+\lfloor un\rfloor}}{x_{m+n}}
			=\left(\frac{m+\lfloor un\rfloor}{m+n}\right)^{-1/\alpha}
			(1+o(1))\longrightarrow u^{-1/\alpha},
		\end{align*}
		it follows that
        \begin{align*}
        x_{m+\lfloor un\rfloor}=u^{-1/\alpha}x_{m+n}(1+o(1)).
		\end{align*}
        Hence for any fixed $u\ge1$, \eqref{h} gives
        \begin{align*}
        \frac{\big|\mu([0,x_{m+\lfloor un\rfloor}])
        -\mu([0,u^{-1/\alpha}x_{m+n}])\big|}
        {\mu([0,x_{m+n}])}
        &\lesssim
        \frac{\left|\int_{u^{-1/\alpha}x_{m+n}}^{x_{m+\lfloor un\rfloor}}x^{-\beta}\,dx\right|}{\int_{0}^{x_{m+n}}x^{-\beta}\,dx}
        \\
        &\lesssim
        \left|
        \left(\frac{x_{m+\lfloor un\rfloor}}{x_{m+n}}\right)^{1-\beta}
        -(u^{-1/\alpha})^{1-\beta}\right|\longrightarrow 0.
		\end{align*}
        Combining this bound with (\ref{mu_rv}), we obtain
		\begin{align*}
			\mu_n^{H}(U_n>u\mid[0,x_{m+n}])
			&=\frac{\mu([0,x_{m+\lfloor un\rfloor}])}{\mu([0,x_{m+n}])}\sim \frac{\mu([0,u^{-1/\alpha}x_{m+n}])}{\mu([0,x_{m+n}])}
			\longrightarrow(u^{-1/\alpha})^{1-\beta}.
		\end{align*}For $u<1$, $\{U_n\leq u\}=\emptyset$.
		Since the distribution of $U$ is continuous, we conclude the proof.
	\end{proof}

\subsubsection*{Step \ref{3}: approximation for $n^{\gamma/\alpha-1}\bigl(S_nf(x)-nf(0)\bigr)$}	
		We first work under
		$\mu_n^{H}(\,\cdot\mid[0,x_{m+n}])$. Write
		\begin{align*}
			Z_n(x):=n^{\gamma/\alpha-1}\bigl(S_nf(x)-nf(0)\bigr).
		\end{align*}
		Fix $\zeta>0$ and $M>1+\zeta$, define
		\begin{align*}
			K_n(\zeta,M):=
			\left\{x\in(0,x_{m+n}]:
			1+\zeta\leq U_n\leq M\right\}.
		\end{align*}
\begin{lemma}\label{approxz_n}As $n\to\infty$,
    \[\sup_{x\in K_n(\zeta,M)}
			\left|Z_n(x)-a c_\alpha^\gamma\frac1n
			\sum_{j=0}^{n-1}\left(U_n(x)-\frac jn\right)^{-\gamma/\alpha}\right|\longrightarrow 0.\]
\end{lemma}
\begin{proof}
		For any $x\in K_n(\zeta,M)$, $T^{nU_n}x\in(x_m,T_L(x_m)]$. Hence, for
		$0\leq j<n$,
		\begin{align}
		\zeta n\leq nU_n-j\leq Mn,	\quad T^jx\in(x_{m+nU_n-j},x_{m+nU_n-j-1}].
			\label{nUn-j}
		\end{align}
		On the one hand, applying Proposition~\ref{xi_xi-1} to \eqref{nUn-j}, we obtain
		\begin{align*}
			\sup_{\substack{x\in K_n(\zeta,M)\\0\leq j<n}}
			\left|\frac{T^jx}{c_\alpha(nU_n(x)-j)^{-1/\alpha}}-1\right|
			\longrightarrow 0.
		\end{align*}
		Applying the $\gamma$-th power still gives
		\begin{align*}
			\sup_{\substack{x\in K_n(\zeta,M)\\0\leq j<n}}
			\left|\frac{(T^jx)^\gamma}
			{c_\alpha^\gamma(nU_n(x)-j)^{-\gamma/\alpha}}-1\right|
			\longrightarrow0.
		\end{align*}
		On the other hand, as $n\to\infty$, the local expansion of $f$ at $0$ gives
		\begin{align*}
			\sup_{\substack{x\in K_n(\zeta,M)\\0\leq j<n}}
			\left|\frac{f(T^jx)-f(0)-a(T^jx)^\gamma}{(T^jx)^\gamma}\right|
			\longrightarrow0.
		\end{align*}
		Combining these two estimates and  $
			(nU_n-j)^{-\gamma/\alpha}\asymp n^{-\gamma/\alpha},$
		we obtain
		\begin{align*}
			&\sup_{\substack{x\in K_n(\zeta,M)\\0\leq j<n}}\left|f(T^jx)-f(0)-a c_\alpha^\gamma(nU_n(x)-j)^{-\gamma/\alpha}\right|\\
			&\le
			\sup_{\substack{x\in K_n(\zeta,M)\\0\leq j<n}}\Biggl[
			\left|\frac{f(T^jx)-f(0)-a(T^jx)^\gamma}{(T^jx)^\gamma}\right|(T^jx)^\gamma\\
			&\quad +|a|c_\alpha^\gamma(nU_n(x)-j)^{-\gamma/\alpha}
			\left|\frac{(T^jx)^\gamma}{c_\alpha^\gamma(nU_n(x)-j)^{-\gamma/\alpha}}-1\right|
			\Biggr]=o(n^{-\gamma/\alpha}).
		\end{align*}
		Thus, for fixed $\zeta,M$, there is 
		$\varepsilon_n\to0$ such that
		\begin{align*}
			\sup_{\substack{x\in K_n(\zeta,M)\\0\leq j<n}}
			\left|f(T^jx)-f(0)-a c_\alpha^\gamma(nU_n(x)-j)^{-\gamma/\alpha}\right|
			\leq\varepsilon_n n^{-\gamma/\alpha}.
		\end{align*}
		It follows that
		\begin{align*}
			\sup_{\substack{x\in K_n(\zeta,M)\\0\leq j<n}}
			\left|n^{\gamma/\alpha-1}\bigl(f(T^jx)-f(0)\bigr)
			-a c_\alpha^\gamma\frac1n\left(U_n(x)-\frac jn\right)^{-\gamma/\alpha}\right|
			\leq \frac{\varepsilon_n}{n}.
		\end{align*}
		Summing over $0\leq j<n$, we have
		\begin{align*}
			&\sup_{x\in K_n(\zeta,M)}
			\left|Z_n(x)-a c_\alpha^\gamma\frac1n
			\sum_{j=0}^{n-1}\left(U_n(x)-\frac jn\right)^{-\gamma/\alpha}\right|\\
			&\leq\sup_{x\in K_n(\zeta,M)}
			\sum_{j=0}^{n-1}
			\left|n^{\gamma/\alpha-1}\bigl(f(T^jx)-f(0)\bigr)
			-a c_\alpha^\gamma\frac1n\left(U_n(x)-\frac jn\right)^{-\gamma/\alpha}\right|\leq\sum_{j=0}^{n-1}\frac{\varepsilon_n}{n}
			=\varepsilon_n\to 0.
		\end{align*}
        \end{proof}

Let $G_n(U_n):=\frac1n\sum_{j=0}^{n-1}\left(U_n-\frac jn\right)^{-\gamma/\alpha}$.
        \subsubsection*{Step \ref{4}: approximate $G_n(U_n)$ by $G_{\gamma/\alpha}(U_n)$}	
		\begin{lemma}\label{15}
		    \[\sup_{x\in K_n(\zeta,M)}
			\left|G_n(U_n(x))-G_{\gamma/\alpha}(U_n(x))\right|\longrightarrow 0,\]
            \begin{align}
			\sup_{x\in K_n(\zeta,M)}
			\left|Z_n(x)-a c_\alpha^\gamma G_{\gamma/\alpha}(U_n(x))\right|
			\longrightarrow0.
			\label{zn_u}
		\end{align}
		\end{lemma}
        \begin{proof}
            Since $U_n\in[1+\zeta,M]$ on $K_n(\zeta,M)$,
		we have
		\begin{align*}
			&\sup_{x\in K_n(\zeta,M)}
			\left|\frac1n\sum_{j=0}^{n-1}
			\left(U_n(x)-\frac jn\right)^{-\gamma/\alpha}-G_{\gamma/\alpha}(U_n(x))\right|\\
			&=\sup_{x\in K_n(\zeta,M)}
			\left|\sum_{j=0}^{n-1}\int_{j/n}^{(j+1)/n}
			\left[\left(U_n(x)-\frac jn\right)^{-\gamma/\alpha}
			-(U_n(x)-t)^{-\gamma/\alpha}\right]d\mathrm{Leb}(t)\right|\\
			&\leq\sup_{x\in K_n(\zeta,M)}
			\sum_{j=0}^{n-1}\int_{j/n}^{(j+1)/n}
			\left|\left(U_n(x)-\frac jn\right)^{-\gamma/\alpha}
			-(U_n(x)-t)^{-\gamma/\alpha}\right|d\mathrm{Leb}(t).
		\end{align*}
		The mean value theorem gives
		\begin{align*}
			&\sup_{x\in K_n(\zeta,M)}
			\sum_{j=0}^{n-1}\int_{j/n}^{(j+1)/n}
			\left|\left(U_n(x)-\frac jn\right)^{-\gamma/\alpha}
			-(U_n(x)-t)^{-\gamma/\alpha}\right|d\mathrm{Leb}(t)\\
			&\leq |\gamma/\alpha|\max_{\zeta\leq v\leq M}v^{-\gamma/\alpha-1}
			\sum_{j=0}^{n-1}\int_{j/n}^{(j+1)/n}
			\left(t-\frac jn\right)d\mathrm{Leb}(t)\\
			&=|\gamma/\alpha|\max_{\zeta\leq v\leq M}v^{-\gamma/\alpha-1}
			\sum_{j=0}^{n-1}\frac1{2n^2}\longrightarrow0.
		\end{align*}
		(\ref{zn_u}) follows from this and Lemma  \ref{approxz_n}.\end{proof}
				
		\subsubsection*{Step \ref{5}: prove the limit distribution}	
        \begin{proof}[Proof of Theorem \ref{snf}\textup{(a)}]
		Fix $\varepsilon$, $\zeta$ and $M$.
		By \eqref{zn_u}, for all sufficiently
		large $n$,
		\begin{align*}
			\mu_n^{H}\left(
			\left|Z_n(x)-a c_\alpha^\gamma G_{\gamma/\alpha}(U_n(x))\right|>\varepsilon,
			x\in K_n(\zeta,M)\mid[0,x_{m+n}]\right)=0.
		\end{align*}
		Hence
		\begin{align*}
			&\mu_n^{H}\left(
			\left|Z_n-a c_\alpha^\gamma G_{\gamma/\alpha}(U_n)\right|>\varepsilon
			\mid[0,x_{m+n}]\right)\\
			&\quad=\mu_n^{H}\left(
			\left|Z_n(x)-a c_\alpha^\gamma G_{\gamma/\alpha}(U_n(x))\right|>\varepsilon,
			x\notin K_n(\zeta,M)\mid[0,x_{m+n}]\right)\\
			&\quad\leq\mu_n^{H}(U_n\leq1+\zeta\mid[0,x_{m+n}])
			+\mu_n^{H}(U_n>M\mid[0,x_{m+n}]).
		\end{align*}
		Since $U_n\to_d U$ under
		$\mu_n^{H}(\,\cdot\mid[0,x_{m+n}])$, this gives
		\begin{align*}
			\limsup_{n\to\infty}\mu_n^{H}\left(
			\left|Z_n-a c_\alpha^\gamma G_{\gamma/\alpha}(U_n)\right|>\varepsilon
			\mid[0,x_{m+n}]\right)
			&\leq\mathbb P(U\leq1+\zeta)
			+\mathbb P(U>M)\\
			&=1-(1+\zeta)^{-\frac{1-\beta}{\alpha}}+M^{-\frac{1-\beta}{\alpha}}.
		\end{align*}
		Letting $M\to\infty$ and $\zeta\downarrow0$, we have, as $n\to\infty$,
		\begin{align}
			Z_n-a c_\alpha^\gamma G_{\gamma/\alpha}(U_n)\longrightarrow0
			\quad \text{in }\mu_n^{H}(\,\cdot\mid[0,x_{m+n}]).
			\label{zn_p}
		\end{align}
		Since $U_n\longrightarrow_d U$ under
		$\mu_n^{H}(\,\cdot\mid[0,x_{m+n}])$, we have
		\begin{align*}
			a c_\alpha^\gamma G_{\gamma/\alpha}(U_n)
			\longrightarrow_d a c_\alpha^\gamma G_{\gamma/\alpha}(U)
			\quad \text{under }\mu_n^{H}(\,\cdot\mid[0,x_{m+n}]).
		\end{align*}
		Combining this with \eqref{zn_p}, we obtain
		\begin{align*}
			Z_n\longrightarrow_d a c_\alpha^\gamma G_{\gamma/\alpha}(U)
			\quad \text{under }\mu_n^{H}(\,\cdot\mid[0,x_{m+n}]).
		\end{align*}
		We transfers this limit to $\mu_n^{H}$ by applying 
		Lemma \ref{en_tv}:
		\begin{align*}
			&\left|\mu_n^{H}(Z_n\leq t)
			-\mu_n^{H}(Z_n\leq t\mid[0,x_{m+n}])\right|\\
			&\quad\leq\left\|\mu_n^{H}(\cdot)
			-\mu_n^{H}(\,\cdot\mid[0,x_{m+n}])\right\|_{\TV}
			\lesssim_{H,\mu}n^{-\beta/\alpha}+n^{-1}\longrightarrow0.
		\end{align*}
		This completes the proof of Theorem \ref{snf}\textup{(a)}.
                \end{proof}

	\subsubsection{Proof of Theorem \ref{snf}\textup{(b)}}
	When $\beta=0$, according to Remark \ref{1and2overlap}, the second class may contribute to the survival probability of the first class at the same order $n^{-1/\alpha}$. The main target of this subsection is to prove that, under \eqref{jmu_rv}, $U_n\longrightarrow_d U$ under $\mu_n^{H}$, the same conclusion as Lemma \ref{un}. The remainder of the proof follows the same strategy as for $0<\beta<1$. 

    We now start with two technical lemmas.

	\begin{lemma}
		\label{jmu_t}
		 $J_\mu(t)\asymp t$ as
		$t\downarrow0$.
	\end{lemma}
\begin{proof}
		For $0<t\leq x_m$,
		the definition of $h_s^{H}$ gives
		\begin{align*}
			\mu\left(E_s^{H}\cap
			\left\{x:\frac12<T^sx\leq\frac{1+t}{2}\right\}\right)
			=\int_{1/2}^{(1+t)/2}h_s^{H}(y)\,d\mathrm{Leb}(y)
			=\int_0^t e_{s+1}^{H}(x)\,d\mathrm{Leb}(x).
		\end{align*}
		Hence we have
		\begin{align*}
			J_\mu(t)=\sum_{s\geq0}\int_0^t e_s^{H}(x)\,d\mathrm{Leb}(x).
		\end{align*}
		Besides, \eqref{h} with $\beta=0$
		gives $t\lesssim_\mu \mu([0,t])\lesssim_\mu t$ for small $t$.
		Thus,
		\begin{align*}
			t
			\lesssim_\mu \mu([0,t])\leq J_\mu(t)\leq t\sum_{s\geq0}\|e_s^{H}\|_{L^\infty([0,x_m])}
			\lesssim_{H,\mu} t,
		\end{align*}
		where we use
		\begin{align*}
			\sum_{s\geq0}\|e_s^{H}\|_{L^\infty([0,x_m])}
		\lesssim_{H,\mu} 1+\sum_{s\geq1}s^{-\frac{1}{\alpha}-1}<\infty.
		\end{align*}
		So $J_\mu(t)\asymp t$ near $0$.

	\end{proof}

	\begin{lemma}
		\label{en_jmu}
For any $u \ge 1$,	
\begin{align}			\mu(E_n^{H})&=J_\mu(x_{m+n})+o(x_{m+n}),\label{11}\\
			\mu_n^{H}(U_n>u)
			&=\frac{J_\mu(x_{m+\lfloor un\rfloor})}{J_\mu(x_{m+n})}+o(1).\label{12}
		\end{align}
	\end{lemma}

	\begin{proof}
		By (\ref{es}), 
		\begin{align}\label{13}
			\mu(E_n^{H})
			=\int_Yh_n^{H}(x)\,d\mathrm{Leb}(x)
			+\sum_{s=0}^n\int_0^{x_{m+n-s}}e_s^{H}(x)\,d\mathrm{Leb}(x),
		\end{align}
		and by Proposition~\ref{en}
		\begin{align*}
			\int_Yh_n^{H}\,d\mathrm{Leb}(x)\lesssim_{H,\mu} n^{-\frac{1}{\alpha}-1}=o(x_{m+n}),
		\end{align*}

Now we prove, for any fixed $u \ge 1$, \begin{align}\label{approxju}
			\sum_{s=0}^n\int_0^{x_{m+\lfloor un\rfloor-s}}e_s^{H}(x)\,d\mathrm{Leb}(x)-J_\mu(x_{m+\lfloor un\rfloor})=o(x_{m+n}).
		\end{align}
        
		Recalling that
		\begin{align*}
			J_\mu(x_{m+\lfloor un\rfloor})=\sum_{s\geq0}\int_0^{x_{m+\lfloor un\rfloor}} e_s^{H}(x)\,d\mathrm{Leb}(x),
		\end{align*}
		so our goal becomes
		\begin{align*}
			\frac1{x_{m+n}}\left|
			\sum_{s=0}^n\int_0^{x_{m+\lfloor un\rfloor-s}}e_s^{H}(x)\,d\mathrm{Leb}(x)
			-\sum_{s=0}^\infty\int_0^{x_{m+\lfloor un\rfloor}}e_s^{H}(x)\,d\mathrm{Leb}(x)\right|
			=o(1).
		\end{align*}
		
		For an integer $0\leq K<n$,
		\begin{align*}
			&\frac1{x_{m+n}}\left|
			\sum_{s=0}^n\int_0^{x_{m+\lfloor un\rfloor-s}}e_s^{H}(x)\,d\mathrm{Leb}(x)
			-\sum_{s=0}^\infty\int_0^{x_{m+\lfloor un\rfloor}}e_s^{H}(x)\,d\mathrm{Leb}(x)\right|\\
			&=\frac1{x_{m+n}}\Bigg|
			\sum_{s=0}^K\int_{x_{m+\lfloor un\rfloor}}^{x_{m+\lfloor un\rfloor-s}}e_s^{H}(x)\,d\mathrm{Leb}(x)+\sum_{s=K+1}^n\int_0^{x_{m+\lfloor un\rfloor-s}}e_s^{H}(x)\,d\mathrm{Leb}(x)\\
			&\qquad-\sum_{s=K+1}^\infty\int_0^{x_{m+\lfloor un\rfloor}}e_s^{H}(x)\,d\mathrm{Leb}(x)
			\Bigg|.
		\end{align*}
		
        Using Proposition~\ref{hn_es} and \eqref{es_k} with $\beta=0$, we continue our estimate
		\begin{align*}
			&\lesssim \left(\sup_{s\geq0}
			\|e_s^{H}\|_{L^\infty([0,x_m])}\right)
			\sum_{s=0}^K\left|\frac{x_{m+\lfloor un\rfloor-s}}{x_{m+\lfloor un\rfloor}}-1\right|\\
            &\qquad +\frac1{x_{m+n}}
			\sum_{s=K+1}^{\lfloor un\rfloor}\int_0^{x_{m+\lfloor un\rfloor-s}}e_s^{H}(x)\,d\mathrm{Leb}(x)+\sum_{s=K+1}^\infty
			\|e_s^{H}\|_{L^\infty([0,x_m])}\\
			&\lesssim_u \sum_{s=0}^K
			\left|\frac{x_{m+n-s}}{x_{m+n}}-1\right|
			+(K+1)^{-\frac1\alpha}+n^{-1}.
		\end{align*}
		
        Since $\frac{x_{m+n-s}}{x_{m+n}}\to1$ as $n\to\infty$.
		Let $n\to\infty$ then let $K\to\infty$, we obtain
		\begin{align*}
			\limsup_{n\to\infty}\frac1{x_{m+n}}\left|
			\sum_{s=0}^n\int_0^{x_{m+\lfloor un\rfloor-s}}e_s^{H}(x)\,d\mathrm{Leb}(x)
			-J_\mu(x_{m+\lfloor un\rfloor})\right|
			\lesssim_u (K+1)^{-\frac1\alpha}\to 0.
		\end{align*}
		
        Thus, we proved (\ref{approxju}). (\ref{approxju}) with $u=1$ and (\ref{13}) prove (\ref{11}).

		Now we turn to prove (\ref{12}). By Lemma \ref{jmu_t},
		\begin{align*}
			\mu(E_n^{H})\asymp J_\mu(x_{m+n})\asymp x_{m+n}.
		\end{align*}

Besides, for any $s\ge 1$,
\begin{align*}
			\mu(E_n^{H}\cap\{U_n> u,\ \sigma_n=0\})
			&=\mu([0,x_{m+\lfloor nu \rfloor}])\\
			&=\int_0^{x_{m+\lfloor nu \rfloor}}h(x)\,d\mathrm{Leb}(x)
			=\int_0^{x_{m+\lfloor nu \rfloor}}e_0^{H}(x)\,d\mathrm{Leb}(x).
		\end{align*}\begin{align*}
			\mu(E_n^{H}\cap\{U_n> u,\ \sigma_n=s\})
			&=\mu\left(E_{s-1}^{H}\cap
			\left\{\frac12<T^{s-1}x\leq\frac{1+x_{m+\lfloor nu \rfloor-s}}2\right\}\right)\\
			&=\int_{1/2}^{(1+x_{m+\lfloor nu \rfloor-s})/2}h_{s-1}^{H}(x)\,d\mathrm{Leb}(x)\\
			&=\int_0^{x_{m+\lfloor nu \rfloor-s}}e_s^{H}(x)\,d\mathrm{Leb}(x).
		\end{align*}

Thus\begin{align*}
			\mu(E_n^{H}\cap\{U_n>u\})
		=\sum_{s=0}^n\int_0^{x_{m+\lfloor un\rfloor-s}}
		e_s^{H}(x)\,d\mathrm{Leb}(x).
		\end{align*}
        
        By Lemma \ref{jmu_t} and $\mu(E_n^H)\asymp x_{m+n}$, we have
		\begin{align*}
			&\left|\mu_n^{H}(U_n>u)
			-\frac{J_\mu(x_{m+\lfloor un\rfloor})}{J_\mu(x_{m+n})}\right|
			=\left|\frac{\mu(E_n^{H}\cap\{U_n>u\})}{\mu(E_n^{H})}
			-\frac{J_\mu(x_{m+\lfloor un\rfloor})}{J_\mu(x_{m+n})}\right|\\
			&\leq\frac{\left|\mu(E_n^{H}\cap\{U_n>u\})
				-J_\mu(x_{m+\lfloor un\rfloor})\right|}{\mu(E_n^{H})}
			+\frac{J_\mu(x_{m+\lfloor un\rfloor})}{J_\mu(x_{m+n})}
			\frac{\left|\mu(E_n^{H})-J_\mu(x_{m+n})\right|}{\mu(E_n^{H})}\\
			&\lesssim \frac1{x_{m+n}}\left|
			\sum_{s=0}^n\int_0^{x_{m+\lfloor un\rfloor-s}}
			e_s^{H}(x)\,d\mathrm{Leb}(x)-J_\mu(x_{m+\lfloor un\rfloor})\right|+o(1)=o(1)
		\end{align*}where the last $o(1)$ is due to (\ref{approxju}). Hence we prove (\ref{12}).
	\end{proof}
Now we prove our main target using two lemmas above.
	\begin{lemma}
		\label{un0}
		Assume \eqref{jmu_rv}, then as $n\to\infty$,
		\begin{align*}
			U_n\longrightarrow_d U
			\quad \text{under }\mu_n^{H}.
		\end{align*}
	\end{lemma}

	\begin{proof}
		Fix $u\geq1$. As $n\to\infty$, Proposition~\ref{xi_xi-1} gives
		\begin{align*}
			\frac{x_{m+\lfloor un\rfloor}}{x_{m+n}}
			=\left(\frac{m+\lfloor un\rfloor}{m+n}\right)^{-1/\alpha}
			(1+o(1))\longrightarrow u^{-1/\alpha}.
		\end{align*}
		Recall that
		\begin{align*}
			J_\mu(t)=\sum_{s\geq0}\int_0^t e_s^{H}(x)\,d\mathrm{Leb}(x).
		\end{align*}
		For $0<a,b\leq x_m$, this gives
		\begin{align*}
			|J_\mu(a)-J_\mu(b)|
			\leq |a-b|\sum_{s\geq0}\|e_s^H\|_{L^\infty([0,x_m])}
			\lesssim |a-b|.
		\end{align*}

		Together with Lemma \ref{jmu_t}, this gives, as $n\to\infty$,
		\begin{align*}
			\left|\frac{J_\mu(x_{m+\lfloor un\rfloor})}
			{J_\mu(x_{m+n})}
			-\frac{J_\mu(u^{-1/\alpha}x_{m+n})}{J_\mu(x_{m+n})}\right|
			\lesssim \left|\frac{x_{m+\lfloor un\rfloor}}{x_{m+n}}-u^{-1/\alpha}\right|
			\longrightarrow0.
		\end{align*}
        
		Thus, using Lemma \ref{en_jmu} and \eqref{jmu_rv}, we have
		\begin{align*}
			\mu_n^{H}(U_n>u)
			=\frac{J_\mu(x_{m+\lfloor un\rfloor})}
			{J_\mu(x_{m+n})}+o(1)=\frac{J_\mu(u^{-1/\alpha}x_{m+n})}{J_\mu(x_{m+n})}+o(1)
			\longrightarrow u^{-1/\alpha}.
		\end{align*} For $u<1$, $\{U_n\leq u\}=\emptyset$. Since $U$ is continuous, this completes the proof.
	\end{proof}

	We now prove Theorem \ref{snf}\textup{(b)} using the same strategy as Theorem \ref{snf}\textup{(a)}.
	\begin{proof}[Proof of Theorem \ref{snf}\textup{(b)}]
		Still let
		\begin{align*}
			Z_n(x):=n^{\gamma/\alpha-1}(S_nf(x)-nf(0)).
		\end{align*}
		Fix $K\geq0$, $\zeta>0$, and $M>1+\zeta$. For $K<n$, define
		\begin{align*}
			K_n(K,\zeta,M):=
			\left\{x\in(0,1]:
			T^nx\in[0,x_m],\ \sigma_n(x)\leq K,\
			1+\zeta\leq U_n(x)\leq M\right\}.
		\end{align*}
		For any $x\in K_n(K,\zeta,M)$, $T^{nU_n}x\in(x_m,T_L(x_m)]$. And for any 
		$\sigma_n\leq j<n$,
		\begin{align*}
		\zeta n\leq nU_n-j\leq Mn,	\quad T^jx\in(x_{m+nU_n-j},x_{m+nU_n-j-1}].
		\end{align*} Hence the arguments in the proof of Lemma \ref{approxz_n} still works when $\sigma_n\leq j<n$. Thus, for fixed $K,\zeta,M$, there exists a sequence $\varepsilon_n\geq0$ with $\varepsilon_n\to0$ such that
		\begin{align*}
			\sup_{\substack{x\in K_n(K,\zeta,M)\\\sigma_n\leq j<n}}
			\left|n^{\gamma/\alpha-1}\bigl(f(T^jx)-f(0)\bigr)
			-a c_\alpha^\gamma\frac1n
			\left(U_n(x)-\frac jn\right)^{-\gamma/\alpha}\right|
			\leq\frac{\varepsilon_n}{n}.
		\end{align*}

		When $x\in K_n(K,\zeta,M)$, the proof of Lemma \ref{15} still works, i.e.,
		\begin{align*}
			\sup_{x\in K_n(K,\zeta,M)}
			\left|\frac1n\sum_{j=0}^{n-1}
			\left(U_n(x)-\frac jn\right)^{-\gamma/\alpha}-G_{\gamma/\alpha}(U_n(x))\right|
			\leq \frac{\gamma}{2\alpha n}\zeta^{-\gamma/\alpha-1}\longrightarrow0.
		\end{align*}
		
		Thus, using $\gamma<\alpha$, we have
		\begin{align}
			&\sup_{x\in K_n(K,\zeta,M)}
			\left|Z_n(x)-a c_\alpha^\gamma G_{\gamma/\alpha}(U_n(x))\right|\notag\\
			&\le\sup_{x\in K_n(K,\zeta,M)}
			\Bigg[
            n^{\gamma/\alpha-1}\sum_{j=0}^{\sigma_n-1}|f(T^jx)-f(0)| \label{16}\\
			&\qquad+
			\sum_{j=\sigma_n}^{n-1}
			\left|n^{\gamma/\alpha-1}\bigl(f(T^jx)-f(0)\bigr)
			-a c_\alpha^\gamma\frac1n
			\left(U_n(x)-\frac jn\right)^{-\gamma/\alpha}\right|
			\label{17}\\
			&\qquad+
			\frac{|a|c_\alpha^\gamma}{n}
			\sum_{j=0}^{\sigma_n-1}\left(U_n(x)-\frac jn\right)^{-\gamma/\alpha}+
			|a|c_\alpha^\gamma
			\left|\frac1n\sum_{j=0}^{n-1}
			\left(U_n(x)-\frac jn\right)^{-\gamma/\alpha}-G_{\gamma/\alpha}(U_n(x))\right|
			\Bigg]\notag\\
			&\leq
			2\|f\|_\infty K n^{\gamma/\alpha-1}
			+\sum_{j=0}^{n-1}\frac{\varepsilon_n}{n}+|a|c_\alpha^\gamma\frac Kn\zeta^{-\gamma/\alpha}
			+\frac{|a|c_\alpha^\gamma\gamma}{2\alpha n}\zeta^{-\gamma/\alpha-1}\longrightarrow0.\notag
		\end{align}

		Fix $\varepsilon>0$, for large $n$, we have
		\begin{align*}
			\mu_n^{H}\left(
		|Z_n(x)-a c_\alpha^\gamma G_{\gamma/\alpha}(U_n(x))|>\varepsilon,\
		x\in K_n(K,\zeta,M)\right)=0.
		\end{align*}
        
Hence for large $n$, \begin{align*}
			\mu_n^{H}\left(
			|Z_n-a c_\alpha^\gamma G_{\gamma/\alpha}(U_n)|>\varepsilon\right)
			&=
			\mu_n^{H}\left(
			|Z_n(x)-a c_\alpha^\gamma G_{\gamma/\alpha}(U_n(x))|>\varepsilon,\
			x\notin K_n(K,\zeta,M)\right)\\
			&\leq
			\mu_n^{H}(T^nx\in Y)+\mu_n^{H}(\sigma_n>K)\\
			&\qquad+\mu_n^{H}(U_n\leq1+\zeta)+\mu_n^{H}(U_n>M).
		\end{align*}
		 
         Therefore Proposition \ref{sigma_un} gives
		\begin{align*}
			\mu_n^{H}\left(
			|Z_n-a c_\alpha^\gamma G_{\gamma/\alpha}(U_n)|>\varepsilon\right)
			&\lesssim \bigl((K+1)^{-1/\alpha}+n^{-1}\bigr)\\
			&\qquad+\mu_n^{H}(U_n\leq1+\zeta)+\mu_n^{H}(U_n>M).
		\end{align*}
        
		Then, by Lemma \ref{un0}, we have
		\begin{align*}
			\limsup_{n\to\infty}
		\mu_n^{H}\left(
		|Z_n-a c_\alpha^\gamma G_{\gamma/\alpha}(U_n)|>\varepsilon\right)
		\leq
		C_{H,\mu}(K+1)^{-1/\alpha}
		+1-(1+\zeta)^{-1/\alpha}+M^{-1/\alpha}.
		\end{align*}

		Letting $K\to\infty$, $M\to\infty$, and
		$\zeta\downarrow0$,
		\begin{align*}
			Z_n-a c_\alpha^\gamma G_{\gamma/\alpha}(U_n)\longrightarrow0
			\quad\text{in }\mu_n^{H}.
		\end{align*}
        
		By Lemma \ref{un0} again, we have
		\begin{align*}
			Z_n\longrightarrow_d a c_\alpha^\gamma G_{\gamma/\alpha}(U)
			\quad \text{under }\mu_n^{H}.
		\end{align*}
	\end{proof}

    \subsection{Proof of Theorem \ref{critical_boundedness}}
We give a complement of Theorem \ref{snf}\textup{(b)}, i.e., consider $\beta=0$ but $\gamma \ge \alpha$.

\begin{proof}[Proof of Theorem \ref{critical_boundedness}]
We divide the integral into three parts:
	\begin{align*}
		\mu_n^H\bigl(S_n|f-f(0)|\bigr)
		&=\int \sum_{j=0}^{n-1}|f(T^jx)-f(0)|\,d\mu_n^H(x)\\
		&=\int_{\{T^nx\in Y\}}\sum_{j=0}^{n-1}|f(T^jx)-f(0)|\,d\mu_n^H(x)\\
		&\qquad +\int_{\{T^nx\in[0,x_m]\}}\sum_{j=0}^{\sigma_n-1}|f(T^jx)-f(0)|\,d\mu_n^H(x)\\
		&\qquad +\int_{\{T^nx\in[0,x_m]\}}\sum_{j=\sigma_n}^{n-1}|f(T^jx)-f(0)|\,d\mu_n^H(x)\\
		&=:I_1+I_2+I_3.
	\end{align*}
	
	We first estimate \(I_1\). Using \(|f-f(0)|\le2\|f\|_\infty\),
	\begin{align*}
		I_1 \le 2n\|f\|_\infty\int_{\{T^nx\in Y\}}\,d\mu_n^H(x)
		=2\|f\|_\infty\,n\,\mu_n^H(T^nx\in Y)\lesssim 1,
	\end{align*}
	where the last inequality follows from Proposition \ref{sigma_un} with \(\beta=0\).

	Next, for \(I_2\),
	\begin{align*}
		I_2
		&\le 2\|f\|_\infty\int_{\{T^nx\in[0,x_m]\}}\sigma_n\,d\mu_n^H(x)
		\le 2\|f\|_\infty\,\mu_n^H(\sigma_n)=2\|f\|_\infty\sum_{K=0}^{n-1}\mu_n^H(\sigma_n>K)\lesssim 1,
	\end{align*}
	where the last inequality also follows from Proposition \ref{sigma_un} with \(\beta=0\).
	
	Finally, we estimate \(I_3\). Since
	\begin{align*}
		\{T^nx\in[0,x_m]\}
		=\bigcup_{s=0}^{n}\{T^nx\in[0,x_m],\,\sigma_n=s\},
	\end{align*}
	\eqref{es} gives
	\begin{align*}
		I_3
		&=\sum_{s=0}^{n}\int_{\{T^nx\in[0,x_m],\,\sigma_n=s\}}
		\sum_{j=s}^{n-1}|f(T^jx)-f(0)|\,d\mu_n^H(x)\\
		&=\frac{1}{\mu(E_n^H)}\sum_{s=0}^{n}\int_0^{x_{m+n-s}}
		\sum_{j=0}^{n-s-1}|f(T_L^jx)-f(0)|\,e_s^H(x)\,d\mathrm{Leb}(x).
	\end{align*}
	By Proposition \ref{hn_es}, we have \(\|e_s^H\|_{L^\infty([0,x_m])}\lesssim (s+1)^{-1/\alpha-1}\) for all \(s\ge0\), and hence
	\begin{align*}
		I_3
		\lesssim \frac{1}{\mu(E_n^H)}\sum_{s=0}^{n}(s+1)^{-1/\alpha-1}
		\int_0^{x_{m+n-s}}\sum_{j=0}^{n-s-1}|f(T_L^jx)-f(0)|\,d\mathrm{Leb}(x).
	\end{align*}
	
	We now consider the inner integral. For \(x\in(x_{m+r+1},x_{m+r}]\) and \(0\le j<n-s\), we have \(T_L^jx\in(x_{m+r-j+1},x_{m+r-j}]\). By  Proposition \ref{xi_xi-1},
	\begin{align*}
		T_L^jx\lesssim (r-j+1)^{-1/\alpha}.
	\end{align*}
	Moreover, \(|f(x)-f(0)|\lesssim x^\gamma\) gives
	\begin{align*}
		|f(T_L^jx)-f(0)|\lesssim (T_L^jx)^\gamma \lesssim (r-j+1)^{-\gamma/\alpha}.
	\end{align*}
    
Using \((0,x_{m+n-s}]=\bigcup_{r \ge n-s}(x_{m+r+1},x_{m+r}]\) and $|x_{m+r+1}-x_{m+r}|\lesssim (r+1)^{-1/\alpha-1}$, we have
\begin{align*}
	\int_0^{x_{m+n-s}}\sum_{j=0}^{n-s-1}|f(T_L^jx)-f(0)|\,d\mathrm{Leb}(x)
	&=\sum_{r=n-s}^{\infty}\int_{x_{m+r+1}}^{x_{m+r}}\sum_{j=0}^{n-s-1}|f(T_L^jx)-f(0)|\,d\mathrm{Leb}(x)\\
	&\lesssim \sum_{r=n-s}^{\infty}(r+1)^{-1/\alpha-1}\sum_{j=0}^{n-s-1}(r-j+1)^{-\gamma/\alpha}.
\end{align*}
Splitting $\sum_{r=n-s}^{\infty}$ at \(r=2(n-s)\), we continue our estimate
	\begin{align*}
		&\lesssim\left(\sum_{r=n-s}^{2(n-s)-1}+\sum_{r=2(n-s)}^{\infty}\right)(r+1)^{-1/\alpha-1}\sum_{j=0}^{n-s-1}(r-j+1)^{-\gamma/\alpha}\\
		&\le (n-s+1)^{-1/\alpha-1}\sum_{r=n-s}^{2(n-s)-1}\sum_{j=0}^{n-s-1}\frac{1}{r-j+1}
		+\sum_{r=2(n-s)}^{\infty}(r+1)^{-1/\alpha-1}\sum_{j=0}^{n-s-1}\frac{1}{n-s}\\
		&\lesssim (n-s+1)^{-1/\alpha-1}(n-s)+(n-s+1)^{-1/\alpha}\\
		&\lesssim (n-s+1)^{-1/\alpha}.
	\end{align*}

	Substituting this bound into \(I_3\) and using \(\mu(E_n^H) \geq\mu([0,x_{m+n}]) \gtrsim (n+1)^{-1/\alpha}\), we have
	\begin{align*}
		I_3
		&\lesssim \frac{1}{\mu(E_n^H)}\sum_{s=0}^{n}(s+1)^{-1/\alpha-1}(n-s+1)^{-1/\alpha}\\
		&\lesssim (n+1)^{1/\alpha}\left(\sum_{0\le s\le n/2}+\sum_{n/2<s\le n}\right)(s+1)^{-1/\alpha-1}(n-s+1)^{-1/\alpha}\lesssim 1.
	\end{align*}
	Combining the estimates for \(I_1\), \(I_2\), and \(I_3\) completes the proof.
\end{proof}

	\subsection{Proof of Theorem \ref{snf_s}}
	Our proof is closely identical to the proof of Theorem \ref{snf}.

	\subsubsection{Proof of Theorem \ref{snf_s}\textup{(a)}}

	\begin{proof}[Proof of Theorem \ref{snf_s}\textup{(a)}] Our proof literally  follows the proof of Theorem \ref{snf}\textup{(a)}, except Lemma \ref{approxz_n} which uses a different local expansion of $f$. Here is a modified proof.

    Let $Z_n(x):=n^{\gamma/\alpha-1}S_nf(x)$. For any $x\in K_n(\zeta,M)$, (\ref{nUn-j}) and \begin{align*}
			\sup_{\substack{x\in K_n(\zeta,M)\\0\leq j<n}}
			\left|\frac{(T^jx)^\gamma}
			{c_\alpha^\gamma(nU_n(x)-j)^{-\gamma/\alpha}}-1\right|
			\longrightarrow0
		\end{align*} still hold. As $n\to\infty$, the local expansion of $f$ near $0$ gives
		\begin{align*}
			\sup_{\substack{x\in K_n(\zeta,M)\\0\leq j<n}}
			\left|\frac{f(T^jx)-a(T^jx)^\gamma}{(T^jx)^\gamma}\right|
			\longrightarrow0.
		\end{align*}
        
		Using $
			(nU_n-j)^{-\gamma/\alpha}\asymp n^{-\gamma/\alpha},$
		we obtain
		\begin{align*}
			&\sup_{\substack{x\in K_n(\zeta,M)\\0\leq j<n}}\left|f(T^jx)-a c_\alpha^\gamma(nU_n(x)-j)^{-\gamma/\alpha}\right|\\
			&\le
			\sup_{\substack{x\in K_n(\zeta,M)\\0\leq j<n}}\Biggl[
			\left|\frac{f(T^jx)-a(T^jx)^\gamma}{(T^jx)^\gamma}\right|(T^jx)^\gamma\\
			&\quad +|a|c_\alpha^\gamma(nU_n(x)-j)^{-\gamma/\alpha}
			\left|\frac{(T^jx)^\gamma}{c_\alpha^\gamma(nU_n(x)-j)^{-\gamma/\alpha}}-1\right|
			\Biggr]=o(n^{-\gamma/\alpha}).
		\end{align*}
        
		Thus, for fixed $\zeta,M$, there is 
		$\varepsilon_n\to0$ such that
		\begin{align*}
			\sup_{\substack{x\in K_n(\zeta,M)\\0\leq j<n}}
			\left|f(T^jx)-a c_\alpha^\gamma(nU_n(x)-j)^{-\gamma/\alpha}\right|
			\leq\varepsilon_n n^{-\gamma/\alpha},
		\end{align*}which implies that
		\begin{align*}
			\sup_{\substack{x\in K_n(\zeta,M)\\0\leq j<n}}
			\left|n^{\gamma/\alpha-1} f(T^jx)
			-a c_\alpha^\gamma\frac1n\left(U_n(x)-\frac jn\right)^{-\gamma/\alpha}\right|
			\leq \frac{\varepsilon_n}{n}.
		\end{align*}
        
		Summing over $0\leq j<n$, we have
		\begin{align*}
			&\sup_{x\in K_n(\zeta,M)}
			\left|Z_n(x)-a c_\alpha^\gamma\frac1n
			\sum_{j=0}^{n-1}\left(U_n(x)-\frac jn\right)^{-\gamma/\alpha}\right|\\
			&\leq\sup_{x\in K_n(\zeta,M)}
			\sum_{j=0}^{n-1}
			\left|n^{\gamma/\alpha-1}f(T^jx)
			-a c_\alpha^\gamma\frac1n\left(U_n(x)-\frac jn\right)^{-\gamma/\alpha}\right|\leq\sum_{j=0}^{n-1}\frac{\varepsilon_n}{n}
			=\varepsilon_n\to 0.
		\end{align*}
	\end{proof}

	\subsubsection{Proof of Theorem~\ref{snf_s}\textup{(b)}}

	\begin{proof}[Proof of Theorem \ref{snf_s}\textup{(b)}]
		Our proof literally follows the proof of Theorem \ref{snf}\textup{(b)} except that (\ref{16}) and (\ref{17}) need modification, where they contribute to $\sum_{j=0}^{\sigma_n-1}$ and $\sum_{j=\sigma_n}^{n-1}$ respectively. Here is a modified proof.
    
	If $x\in K_n(K,\zeta,M)$ and $0\leq j<\sigma_n$, then $T^jx>x_{m+K}$. Indeed, if $T^jx\leq x_{m+K}$, then $\sigma_n-1-j\leq K-1$ and
		\begin{align*}
			T^{\sigma_n-1}x=T_L^{\sigma_n-1-j}(T^jx)
			\leq x_{m+K-(\sigma_n-1-j)}\leq x_m,
		\end{align*}
	 contradicting the definition of $\sigma_n$. Since $f$ is bounded on $[x_{m+K},1]$, then
		\begin{align*}
			B_K:=\sup_{x_{m+K}\leq x\leq1}|f(x)|<\infty.
		\end{align*}
		Thus, for fixed $K$, as $n\to\infty$,
		\begin{align*}
			\sup_{x\in K_n(K,\zeta,M)}
			n^{\gamma/\alpha-1}
			\sum_{j=0}^{\sigma_n-1}|f(T^jx)|
			\leq KB_Kn^{\gamma/\alpha-1}\longrightarrow0.
		\end{align*}

		If $x\in K_n(K,\zeta,M)$ and $\sigma_n\leq j<n$, we replace $0\le j<n$ by $\sigma_n\leq j<n$ in the modified proof of Theorem \ref{snf_s}\textup{(a)}, giving \begin{align*}
			\sup_{\substack{x\in K_n(K,\zeta,M)\\ \sigma_n\leq j<n}}
			\left|n^{\gamma/\alpha-1} f(T^jx)
			-a c_\alpha^\gamma\frac1n\left(U_n(x)-\frac jn\right)^{-\gamma/\alpha}\right|
			\leq \frac{\varepsilon_n}{n}.
		\end{align*}
	\end{proof}

	\subsection{Proof of Corollary \ref{snf_bv}}

	\begin{proof}[Proof of Corollary \ref{snf_bv}]
		Let $\mathcal P$ be the transfer operator of $T$ with respect to Lebesgue measure, and let $\mathring{\mathcal P}:=\1_{H^c}\mathcal P$. Since $H,H^c$ are finite union of intervals, $E^H_n$ is too for any $n\ge 0$. We claim $h_n^H(0.5+)$ exists for any $n\ge 0$. This is because $h\in \BV([\epsilon,1])$, $h(T^{-n}(0.5+))$ exists, hence \[h_{n}^H(0.5+)=\sum_{T^ny=0.5+}\frac{\1_{E^{H}_n}(y)h(y)}{(T^n)'(y)}\]  exists. On the other hand, by Proposition \ref{hn_es},
		\begin{align*}
			\sum_{n\geq0}\|h_n^H\|_{L^\infty(Y)}<\infty.
		\end{align*}
		Hence as $t\downarrow0$,
		\begin{align*}
			\frac{J_\mu(t)-\mu([0,t])}t
			=\sum_{n\geq0}\frac1t\int_{1/2}^{(1+t)/2}h_n^H(y)\,d\mathrm{Leb}(y)
			\longrightarrow \frac12\sum_{n\geq0}h_n^H(1/2+)=:A<\infty,
		\end{align*} i.e., $J_\mu(t)=\mu([0,t])+At+o(t)$. Using $t \lesssim \mu([0,t])\le J_\mu(t)$ and (\ref{mu_rv}), we have

\[
\frac{J_\mu(\lambda t)}{J_\mu(t)}-\lambda
=
\frac{\mu([0,t])}{J_\mu(t)}
\left(
\frac{\mu([0,\lambda t])}{\mu([0,t])}-\lambda
\right)
+\frac{o(t)}{J_\mu(t)}\longrightarrow 0.
\]

This corollary follows from Theorems \ref{snf}\textup{(b)} and \ref{snf_s}\textup{(b)} respectively.
	\end{proof}

\section{Appendix}\label{appendix}
\subsection{\texorpdfstring{$G_{\gamma/\alpha}(U)$}{G[gamma/alpha](U)} is not stable}
By \cite[Definition~1.1.6]{ST94}, the characteristic function of a stable/Gaussian law has the form
\[
 \varphi(t)=
 \begin{cases}
 \displaystyle
 \exp\!\left\{-\sigma^{\vartheta}|t|^{\vartheta}
 \left(1-\mathrm i\eta\operatorname{sgn}(t)
 \tan\frac{\pi\vartheta}{2}\right)+\mathrm i\zeta t\right\},
 &\vartheta\ne1,\\[8pt]
 \displaystyle
 \exp\!\left\{-\sigma|t|
 \left(1+\mathrm i\eta\frac{2}{\pi}\operatorname{sgn}(t)
 \log|t|\right)+\mathrm i\zeta t\right\},
 &\vartheta=1,
 \end{cases}
\]
for some $0<\vartheta\le2$, $\sigma\ge0$, $-1\le\eta\le1$, and
$\zeta\in\mathbb R$. We will study the characteristic function of our limit distribution $G_{\gamma/\alpha}(U)$ and show that this is not a stable law.

Recall that $0<\alpha<1$ and $0\leq\beta<1$. Let $\varphi_{G}$ denote the characteristic function of $G_{\gamma/\alpha}(U)$. Since
$U$ has density $\frac{1-\beta}{\alpha}u^{-(1-\beta)/\alpha-1}$ on
$(1,\infty)$, we have
\[
 \varphi_{G}(t)=\mathbb E e^{\mathrm i tG_{\gamma/\alpha}(U)}
 =\frac{1-\beta}{\alpha}\int_1^\infty
 e^{\mathrm i tG_{\gamma/\alpha}(u)}u^{-(1-\beta)/\alpha-1}\,du.
\]

\begin{proposition}\label{nonstable}
For every $\gamma\ne0$, the distribution of $G_{\gamma/\alpha}(U)$ is not a stable law.
\end{proposition}

\begin{proof}
We first investigate the distribution of $G_{\gamma/\alpha}(U)$ near some endpoints. For $u>1$, we have
\[
 G_{\gamma/\alpha}'(u)=-\frac{\gamma}{\alpha}
 \int_0^1(u-v)^{-\gamma/\alpha-1}\,dv
 =u^{-\gamma/\alpha}-(u-1)^{-\gamma/\alpha}.
\]
Hence $G_{\gamma/\alpha}$ is strictly monotone.
For $\gamma>0$, $G_{\gamma/\alpha}$ decreases to zero as $u\to\infty$,
while for $\gamma<0$, $G_{\gamma/\alpha}$ increases from
$G_{\gamma/\alpha}(1)=\int_0^1(1-v)^{-\gamma/\alpha}\,dv
=\alpha/(\alpha-\gamma)$.

For $\gamma>0$, as $u\to\infty$, we have
\[
 G_{\gamma/\alpha}(u)=u^{-\gamma/\alpha}
 \int_0^1(1-v/u)^{-\gamma/\alpha}\,dv
 \sim u^{-\gamma/\alpha}.
\]
The substitution
$x=G_{\gamma/\alpha}(u)$ gives, as $x\downarrow0$,
\begin{align}\label{0endpoint}
 \mathbb P(G_{\gamma/\alpha}(U)\le x)=\mathbb P(U\ge u)=u^{-(1-\beta)/\alpha}
 \sim x^{(1-\beta)/\gamma}.
\end{align}

For $\gamma<0$, one can check that $G_{\gamma/\alpha}'(u)\to1$ as
$u\downarrow1$. Hence, as $u\downarrow1$,
\[
 G_{\gamma/\alpha}(u)-\frac{\alpha}{\alpha-\gamma}\sim u-1.
\]
Then the substitution
$x=G_{\gamma/\alpha}(u)-\frac{\alpha}{\alpha-\gamma}$ gives, as $x\downarrow0$,
\begin{align}\label{alphaendpoint}
 \mathbb P\left(G_{\gamma/\alpha}(U)\le\frac{\alpha}{\alpha-\gamma}+x\right)
 =\mathbb P(U\le u)=1-u^{-(1-\beta)/\alpha}
 \sim\frac{1-\beta}{\alpha}x.
\end{align}

Next, we argue by contradiction, comparing the endpoint behavior of the distribution
of $G_{\gamma/\alpha}(U)$. We suppose that $G_{\gamma/\alpha}(U)$ has a stable law, then
\[
 |\varphi_{G}(t)|=\exp(-\sigma^{\vartheta}|t|^{\vartheta}),
\]
where $\vartheta\in(0,2]$ and $\sigma>0$ since $G_{\gamma/\alpha}(U)$ is nondegenerate. It follows that for every $k\in\mathbb N_0$,
\[
 \int_{\mathbb R}|t|^k|\varphi_{G}(t)|\,dt<\infty.
\]
Hence by Fourier inversion, $G_{\gamma/\alpha}(U)$ has a density $h_{G}\in C^\infty(\mathbb R)$
given by
\[
 h_{G}(x)=\frac1{2\pi}\int_{\mathbb R}
 e^{-\mathrm i tx}\varphi_{G}(t)\,dt.
\]
If $\gamma>0$, we have $G_{\gamma/\alpha}(U)>0$ almost surely, then the density $h_{G}$ is zero
on $(-\infty,0)$. Similarly, if $\gamma<0$, we have $h_{G}$ is zero on
$(-\infty,\frac{\alpha}{\alpha-\gamma})$. Since $h_{G}\in C^\infty(\mathbb R)$, all derivatives of $h_{G}$ are zero at $0$ when $\gamma>0$,
and at $\frac{\alpha}{\alpha-\gamma}$ when $\gamma<0$. Thus, for every
$N\in\mathbb N$, Taylor's formula gives, as $x\downarrow0$,
\[
 \begin{aligned}
 h_{G}(x)&=O(x^{N-1}), &&\gamma>0,\\
 h_{G}\left(\frac{\alpha}{\alpha-\gamma}+x\right)
 &=O(x^{N-1}), &&\gamma<0.
 \end{aligned}
\]
Integrating from $0$ to $x$ when $\gamma>0$, and from
$\frac{\alpha}{\alpha-\gamma}$ to $\frac{\alpha}{\alpha-\gamma}+x$
when $\gamma<0$, we obtain
\[
 \begin{aligned}
 \mathbb P(G_{\gamma/\alpha}(U)\le x)&=O(x^N), &&\gamma>0,\\
 \mathbb P\left(G_{\gamma/\alpha}(U)\le\frac{\alpha}{\alpha-\gamma}+x\right)
 &=O(x^N), &&\gamma<0.
 \end{aligned}
\]

We can now find contradictions. For $\gamma>0$, choose $N\in\mathbb N$ such that $N>(1-\beta)/\gamma$.
Then as $x\downarrow0$,
\[
 \frac{\mathbb P(G_{\gamma/\alpha}(U)\le x)}{x^{(1-\beta)/\gamma}}
 =O\left(x^{N-(1-\beta)/\gamma}\right)\longrightarrow0,
\]
whereas \eqref{0endpoint} gives a nonzero limit. For $\gamma<0$, taking $N=2$ gives
\[
 \frac{\mathbb P\left(G_{\gamma/\alpha}(U)\le\frac{\alpha}{\alpha-\gamma}+x\right)}{x}
 =O(x)\longrightarrow0,
\]
whereas \eqref{alphaendpoint} gives a nonzero limit. Thus, the
distribution of $G_{\gamma/\alpha}(U)$ is not a stable law, let alone Gaussian law.
\end{proof}

\begin{remark}
In fact, further calculations yield the following asymptotic estimate for the characteristic function. As $|t|\to\infty$,
\[
 \log|\varphi_{G}(t)|=-\kappa\log|t|+O(1),
\]
where
\[
 \kappa:=
 \begin{cases}
 \displaystyle
 \min\!\left\{
 \frac{1-\beta}{\gamma},
 \frac{\alpha}{\alpha-\gamma}
 \right\},
 &0<\gamma<\alpha,\\[3mm]
 \displaystyle\frac{1-\beta}{\gamma},
 &\gamma\ge\alpha,\\[3mm]
 1,
 &\gamma<0.
 \end{cases}
\]
\end{remark}

	\subsection{\texorpdfstring{$n^{\gamma/\alpha-1}$}{n to the power (gamma/alpha - 1)} is a unique normalization order}
	\label{rn}
	In this subsection, we show that $n^{\gamma/\alpha-1}$ is a unique order of normalization in Theorems~\ref{snf} and~\ref{snf_s}.
	Under the assumptions of Theorem \ref{snf}, let $r_n>0$ be any deterministic sequence. Then
	\begin{align}\label{18}
		r_n\bigl(S_nf-nf(0)\bigr)
		=\bigl(r_n n^{1-\gamma/\alpha}\bigr)
		n^{\gamma/\alpha-1}\bigl(S_nf-nf(0)\bigr).
	\end{align}
	The second factor converges in distribution to $Z:=a c_\alpha^\gamma G_{\gamma/\alpha}(U)$. We now consider all possible limit behaviors of $r_n n^{1-\gamma/\alpha}$.
	\begin{enumerate}[label=\textbf{\arabic*.}]
		\item {\bfseries\boldmath $r_n n^{1-\gamma/\alpha}\to0$.} 
        
        For every $\varepsilon>0$, $\frac{\varepsilon}{r_n n^{1-\gamma/\alpha}}\to \infty$. Hence, by Theorem \ref{snf},
	\begin{align*}
			\mu_n^{H}\left(
			\left|r_n\bigl(S_nf-nf(0)\bigr)\right|>\varepsilon
			\right)
			=
			\mu_n^{H}\left(
			\left|n^{\gamma/\alpha-1}\bigl(S_nf-nf(0)\bigr)\right|
			>
			\frac{\varepsilon}{r_n n^{1-\gamma/\alpha}}
			\right)\to 0,
	\end{align*} i.e., as $n\to\infty$,
	\begin{align*}
			r_n\bigl(S_nf-nf(0)\bigr)\longrightarrow0
			\quad\text{in }\mu_n^{H}.
	\end{align*}
		\item {\bfseries\boldmath $r_n n^{1-\gamma/\alpha}\to c\in(0,\infty)$.}

		Since the first factor of (\ref{18}) tends to $c$ and the second converges in distribution to $a c_\alpha^\gamma G_{\gamma/\alpha}(U)$, we obtain as $n\to\infty$,
	\begin{align*}
			r_n\bigl(S_nf-nf(0)\bigr) \longrightarrow_d c a c_\alpha^\gamma G_{\gamma/\alpha}(U)
			\quad\text{under }\mu_n^{H}.
	\end{align*}
		\item {\bfseries\boldmath $r_n n^{1-\gamma/\alpha}\to\infty$.}

		Assume $a>0$. For any $M>0$, $\frac{M}{r_n n^{1-\gamma/\alpha}}\to 0$. Hence, by Theorem \ref{snf},
	\begin{align*}
			\mu_n^{H}\left(
			r_n\bigl(S_nf-nf(0)\bigr)\leq M
			\right)
			=
			\mu_n^{H}\left(
			n^{\gamma/\alpha-1}\bigl(S_nf-nf(0)\bigr)
			\leq
			\frac{M}{r_n n^{1-\gamma/\alpha}}
			\right)\to 0.
	\end{align*}

		Therefore, as $n\to\infty$,
	\begin{align*}
			r_n\bigl(S_nf-nf(0)\bigr)\longrightarrow+\infty
			\quad\text{in }\mu_n^{H}.
	\end{align*}
		If $a<0$, the same argument gives $-\infty$ in $\mu_n^{H}$.

		\item {\bfseries\boldmath $r_n n^{1-\gamma/\alpha}$ does not converge in $[0,\infty]$.}

		If $r_n n^{1-\gamma/\alpha}$ is unbounded, there is a subsequence tending to infinity. By Case 3, the corresponding distributions are not tight on $\mathbb R$.

		If $r_n n^{1-\gamma/\alpha}$ is bounded but does not converge, it has two subsequences with distinct limits $0\leq c_1<c_2<\infty$. By Cases 1--2, the corresponding limit distributions are $c_1Z$ and $c_2Z$, respectively. They are distinct, since $\mathbb P(|Z|>0)=1$. Thus, in either case, the limit distribution does not exist.
	\end{enumerate}

	Thus, $r_n(S_nf-nf(0))$ has a nondegenerate limit distribution under $\mu_n^{H}$ if and only if, as $n\to\infty$,
	\begin{align*}
		r_n n^{1-\gamma/\alpha}\longrightarrow c\in(0,\infty).
	\end{align*}

	For Theorem \ref{snf_s}, the same argument follows literally by replacing   $S_nf-nf(0)$ with $S_nf$.

	\section*{Acknowledgements}
	We thank Prof. Sandro Vaienti for helpful discussions and comments.
	\bibliography{sn-bibliography}
	
\end{document}